\documentclass{amsart}

\usepackage{amsmath, amsthm, amsfonts, ifpdf}
\usepackage{color}
\usepackage{graphicx}
\usepackage{multicol}
\usepackage{verbatim}
\usepackage{tikz}
\usepackage{amssymb}

\usetikzlibrary{arrows, shapes, trees, backgrounds}
\usetikzlibrary{intersections}

\theoremstyle{plain}
\newtheorem{theorem}{Theorem}[section]
\newtheorem*{theorem*}{Theorem}

\newtheorem{pro}[theorem]{Proposition}
\newtheorem{Def}[theorem]{Definition}
\newtheorem{lemma}[theorem]{Lemma}

\newtheorem{Pro}[theorem]{Problem}

\newtheorem{cor}[theorem]{Corollary}

\theoremstyle{definition}
\newtheorem*{Def*}{Definition}
\newtheorem{Rem}[theorem]{Remark}
\newtheorem*{Rem*}{Remark}

\newtheorem*{changelog*}{Changelog}

\numberwithin{equation}{section}

\newcommand{\bpo}{\begin{pro}}
    \newcommand{\epo}{\end{pro}}
\newcommand{\be}{\begin{equation}}
    \newcommand{\ene}{\end{equation}}
\newcommand{\br}{\begin{Rem}}
    \newcommand{\er}{\end{Rem}}
\newcommand{\bl}{\begin{lem}}
    \newcommand{\el}{\end{lem}}
\newcommand{\bd}{\begin{Def}}
    \newcommand{\ed}{\end{Def}}
\newcommand{\ben}{\begin{enumerate}}
    \newcommand{\een}{\end{enumerate}}
\newcommand{\bp}{\begin{proof}}
    \newcommand{\ep}{\end{proof}}
\newcommand{\beq}{\begin{equation*}}
    \newcommand{\eeq}{\end{equation*}}
\newcommand{\bear}{\begin{eqnarray*}}
    \newcommand{\eear}{\end{eqnarray*}}
\newcommand{\bt}{\begin{theorem}}
    \newcommand{\et}{\end{theorem}}
\newcommand{\bst}{\begin{split}}
    \newcommand{\est}{\end{split}}

\newcommand{\bal}{\begin{aligned}}
    \newcommand{\eal}{\end{aligned}}

\newcommand{\F}[2]{\frac{#1}{#2}}

\renewcommand{\dim}{\mathrm{dim}}

\newcommand{\R}{\mathbb{R}}

\newcommand{\PLH}{{\mkern-1mu\times\mkern-1mu}}

\newcommand{\rH}{\mathrm{H}}
\newcommand{\dvol}{\mathrm{dvol}}

\newcommand{\rd}{\mathrm{d}}

\newcommand{\rgr}{\mathrm{gr}}
\renewcommand{\div}{\mathrm{div}}
\newcommand{\spt}{\operatorname{supp}}
\newcommand{\BV}{\mathrm{BV}}

\def\XXint#1#2#3{{\setbox0=\hbox{$#1{#2#3}{\int}$}
        \vcenter{\hbox{$#2#3$}}\kern-.5\wd0}}

\makeatletter
\def\@citestyle{\m@th\upshape\mdseries}
\def\citeform#1{{\bfseries#1}}
\def\@cite#1#2{{%
        \@citestyle[\citeform{#1}\if@tempswa, #2\fi]}}
\@ifundefined{cite }{%
    \expandafter\let\csname cite \endcsname\cite
    \edef\cite{\@nx\protect\@xp\@nx\csname cite \endcsname}%
}{}
\makeatother

\begin{document}

    \title[Topology of PMC hypersurfaces]{Prescribed mean curvature hypersurfaces in conformal product manifolds}
    \author{Qiang Gao and Hengyu Zhou}

    \address[H. Zhou]{College of Mathematics and Statistics, Chongqing University, Huxi Campus, Chongqing 401331, P.\,R.\,China}
    \address{Chongqing Key Laboratory of Analytic Mathematics and Applications, Chongqing University, Huxi Campus, Chongqing 401331, P.\,R.\,China}
    \email{zhouhyu@cqu.edu.cn}

    \address[Q. Gao]{School of Mathematical Sciences, Shenzhen University, Shenzhen, Guangdong 518060, P.\,R.\,China}
    \email{gaoqiangks@szu.edu.cn}


    \begin{abstract}
        In this paper, we establish the existence of prescribed mean curvature (PMC) hypersurfaces in conformal product manifolds with (possibly empty) $C^{1,\alpha}$ fixed graphical boundaries under a barrier condition. This result generalizes Gerhardt's work to non-flat conformal backgrounds. As a consequence, we obtain new solutions to the high-dimensional PMC Plateau problem with explicitly characterized topology. Moreover, under a quasi-decreasing condition on the PMC function, we demonstrate that the resulting hypersurfaces are $C^1$ graphs.
    \end{abstract}

    \date{\today}
    \maketitle

    \section{Introduction}

    This paper is motivated by the study of prescribed mean curvature (PMC) hypersurfaces, a problem originating from Yau~\cite{Yau82}. We refer to this as the Prescribed Mean Curvature (PMC) problem, stated as follows:

    \begin{Pro}[The PMC problem, {\cite[ p.\,683]{Yau82}}]\label{key:problem}
        Let $M^{n+1}$ be an $(n+1)$-dimensional closed or compact Riemannian manifold with $n\geq 2$, and let $\mathcal{H}$ be a $C^{1,\alpha}$ function on its tangent bundle $TM^{n+1}$. Does there exist a hypersurface $\Sigma$ in $M$ whose mean curvature is given by $\mathcal{H}$, and how can one describe its topology?
    \end{Pro}

    In the case where the boundary $\partial \Sigma$ (if nonempty) is closed and fixed, the PMC problem in this setting corresponds to the Plateau problem for PMC hypersurfaces. Here, the prescribed mean curvature function $\mathcal{H} = \mathcal{H}(p, \vec{\nu})$ may depend on both the position $p$ and the unit normal vector $\vec{\nu}$ of the hypersurface.

    Some special cases regarding the existence of solutions to the PMC problem have been extensively studied. When $M^{n+1}$ is closed and $\mathcal{H}=\mathcal{H}(p)$, the existence of minimal hypersurfaces (i.e., $\mathcal{H}\equiv 0$) was established by Pitts~\cite{Pitts81} for $2\leq n\leq 5$, by Schoen and Simon~\cite{SS81} for $n=6$, and by Zhou and Zhu in~\cite{ZZ19} for $\mathcal{H}\equiv C$ and in~\cite{ZZ20} for generic $\mathcal{H}(p)$, using a min-max theory of PMC functionals. For compact manifolds satisfying a barrier condition (see Definition~\ref{definiton:barrier:condition}), Eichmair~\cite{Eich09} established the existence of PMC hypersurfaces with fixed (possibly empty) closed boundaries for general $\mathcal{H}$ via a blow-up method. In the context of positive scalar curvature manifolds, certain such PMC hypersurfaces are also referred to as $\mu$-bubbles~\cite{gro18}.

    Compared to the existence theory for the PMC problem, relatively little is known about the possible topology of solutions. Most results in this direction have been obtained in three-dimensional manifolds for constant $\mathcal{H}$. Relevant work includes the existence of embedded minimal disks by Meeks and Yau~\cite{Yau97} (see also~\cite{HS88}), the existence of PMC disks in sufficiently small geodesic balls (with $\mathcal{H} = \mathcal{H}(p)$) by Gulliver and Spruck~\cite{GS08}, and the existence of embedded constant mean curvature disks under additional assumptions by Coskunuzer~\cite{Cos17}.

    An alternative approach to the PMC problem is based on partial differential equations. Under a barrier assumption, the existence of $C^2$ star-shaped graphs with mean curvature $\mathcal{H}$ satisfying certain monotonicity conditions was shown by Bakelman and Kantor~\cite{BK74} and by Treibergs and Wei~\cite{TW83} in Euclidean spaces, and by Hsu, Shiau, and Wang~\cite{HSW00} in hemispheres. These monotonicity conditions are special cases of \eqref{condition:monotone} in conformal product manifolds. Namely, by Theorem~\ref{at:G:key}, these PMC hypersurfaces in~\cite{BK74, TW83, HSW00} are $C^2$ graphs of the solutions to certain PMC equations on spheres assuming \eqref{condition:monotone}.

    The monotonicity condition \eqref{condition:monotone} plays a crucial role in establishing gradient estimates for solutions to the PMC equation. Nevertheless, such a condition is not always satisfied, even by minimal graphs in conformal product manifolds (see Remark~\ref{nonatural:rk}). In~\cite{Ger96}, Gerhardt removes this monotonicity assumption on the PMC function and proves the existence of closed PMC spheres under barrier conditions in conformally flat spaces via a refined approximation scheme.

    The purpose of this paper is to generalize Gerhardt's work to broader classes of manifolds. We define a conformal product manifold $M_f$ by
    \be \label{def:cpd}
    M_f:=\left\{ N\PLH (a,b), e^{2f}(\sigma+dr^2) \right\},
    \ene
    where $(a,b)$ is an open interval, $N$ is an $n$-dimensional closed or $C^{2,\alpha}$ compact Riemannian manifold with metric $\sigma$, and $f$ is a smooth function defined on a neighborhood of $N\PLH (a,b)$. Examples of conformal product manifolds include Euclidean spaces, hyperbolic spaces, conformally flat spaces, and all warped product manifolds.

    Throughout this paper, given any function $h_1$, we denote by $\mathrm{gr}(h_1)$ the graph of $h_1$. Let $u_1, u_0 \in C^{2}(N^0) \cap C(N)$ be two functions, where $N^0$ denotes the interior of $N$, and let $\mathcal{H}$ be a $C^{1,\alpha}$ function defined on a neighborhood of $TM_f$.

    \begin{Def}\label{def:barrier:basic}
        We say that $(\rgr(u_1), \rgr(u_0))$ forms a barrier for $(M_f, \mathcal{H})$ provided that $a<u_1< u_0 <b$ on $N^0$ and
        \be \label{barrier:condition}
        \rH_{\rgr(u_1)}(p)\leq \mathcal{H}(p,\vec{\nu}),\quad
        \rH_{\rgr(u_0)}(p)\geq  \mathcal{H}(p,\vec{\nu}).
        \ene
        Here $\rH_{\rgr(u_i)}$ is the mean curvature of $\rgr(u_i)$ with respect to its upward unit normal for $i=1,0$. When $N$ is of class $C^{2,\alpha}$ and compact, we require that there exists $\psi\in C^{1,\alpha}(\partial N)$ such that $u_1=u_0=\psi$ on $\partial N$ (see Definition~\ref{definiton:barrier:condition}).
    \end{Def}

    The main result of this paper is stated as follows.

    \begin{theorem}\label{thm:main}
        Suppose $(\rgr(u_1), \rgr(u_0))$ forms a barrier for $(M_f, \mathcal{H})$, and let $2 \leq n \leq 7$. Then there exists a $C^{3,\alpha}$ embedded hypersurface $\Sigma$ with prescribed graphical boundary $\mathrm{gr}(\psi)$ (which may be empty) such that $\Sigma \cup \mathrm{gr}(\psi)$ is homeomorphic to $N$ and has mean curvature $\mathcal{H}$.
    \end{theorem}

    This result generalizes Gerhardt's result~\cite{Ger96} to the case when $N$ has nonempty boundary. The case where $N$ has compact boundary provides explicit examples of solutions to the Plateau problem for PMC hypersurfaces in higher dimensions with well-controlled topology, as posed in Problem~\ref{key:problem}. These examples are nontrivial, since the resulting PMC hypersurfaces may not be $C^1$ graphs over $N$, even as limits of $C^2$ graphs with bounded mean curvature in the $C^{1,\gamma}$ sense. To our knowledge, within the context of the Plateau problem for PMC hypersurfaces, very few examples of this type are known beyond classical solutions to the Dirichlet problem for the PMC equation. The construction of these hypersurfaces (via Theorem~\ref{key:thm}) is also of independent interest in the theory of PMC equations.

    A natural question arising from Theorem~\ref{thm:main} is under what nontrivial conditions the hypersurface $\Sigma$ is a $C^1$ regular graph.
    In Theorem~\ref{quasi:monotone:case}, this question is addressed in part under a quasi-decreasing condition; specifically, by assuming $f \equiv 0$ and that $\mathcal{H}$ admits a decomposition
    $$
    \mathcal{H}(x,r, Y,t)=\mathcal{H}_1(x,r,Y,t)+\mathcal{H}_2(x,r,Y,t)t
    $$
    with the quasi-decreasing condition
    $$
    \F{\partial \mathcal{H}_1}{\partial r}\leq 0
    $$
    for $x\in N$, $r,t\in \R$, $Y\in T_x N$, where $\mathcal{H}_1, \mathcal{H}_2$ are $C^{1,\alpha}$ functions on $T(N\PLH\R)$.

    Our approach to proving Theorem~\ref{thm:main} follows a line of reasoning similar to that in~\cite{Ger96}, where Gerhardt constructed an increasing sequence of $C^2$ open graphs whose limit is a closed embedded hypersurface homeomorphic to $S^n$ with prescribed mean curvature under the dimensional constraint $2 \leq n \leq 6$. The primary difficulty lies in computing the mean curvature of the limit surface $\Sigma$.

    To overcome this difficulty, we introduce an alternative strategy to Gerhardt's integral formula~\cite{Ger96}.
    In Theorem~\ref{key:thm}, we characterize the $L^1$ closure of $C^2$ functions with bounded mean curvature, bounded height, and whose first-order derivatives are bounded in the distributional $L^1$ sense. As an application, we prove that $\Sigma$ is the graph of a function in $W^{1,1}(N)$. Furthermore, we show that over a dense open set $\Omega \subset N$, the restriction $\Sigma|_{\Omega \times \mathbb{R}}$ is a $C^2$ graph with mean curvature $\mathcal{H}$, and is hence dense in $\Sigma$. Combining the regularity of $\Sigma$ with classical elliptic theory, we conclude that the mean curvature of the entire hypersurface $\Sigma$ is $\mathcal{H}$.

    We conclude with several remarks in order. First, this paper is a sequel to~\cite{zhou22a}. The main difference between the two papers is that the latter studies PMC graphs satisfying an Ncf property on their underlying domains, whereas the present article approaches the problem from the perspective of a barrier condition for PMC functions within the geometric setting of conformal product manifolds. Second, the PMC hypersurfaces considered here need not be critical points of area-type functionals. For example, in the case of the product manifold $N \times \mathbb{R}$, one may take $\mathcal{H} = \operatorname{tr}_N(k) - k(\vec{\nu}, \vec{\nu})$, where $k$ is a smooth $(0,2)$-tensor field on $N$. Finally, it would be of considerable interest to resolve Problem~\ref{key:problem} under the barrier condition without assuming conformally product structures.

    This paper is organized as follows. In Section 2, we recall some fundamental facts concerning barrier conditions and $\Lambda$-perimeter minimizers. Section 3 is devoted to the study of limits of $C^2$ graphs with bounded mean curvature and bounded height. In Section 4, we establish the existence of a PMC hypersurface under a barrier condition. Finally, in Section 5, we demonstrate how the quasi-decreasing condition on the prescribed mean curvature function leads to the existence of a $C^1$ PMC graph, thereby completing the proof of Theorem~\ref{thm:main}.

    This paper is supported by the National Natural Science Foundation of China (Grant No.~11801046) and partially sponsored by the Natural Science Foundation of Chongqing, China (Grant No.~cstc2021jcyj-msxmX0430).

    \section{Conformal product manifolds and $\Lambda$-perimeter minimizers}
    This section has two main objectives: first, to introduce a barrier condition for PMC hypersurfaces in conformal product manifolds, and second, to provide preliminary material on $\Lambda$-perimeter minimizers.

    \subsection{Conformal product manifolds}
    Suppose $M$ is an $(n+1)$-dimensional Riemannian manifold with $n\geq 2$, endowed with a metric $g$. Given any smooth function $f$ on $M$, we denote by $M_f$ the conformal manifold $M$ equipped with the conformal metric $e^{2f}g$.

    \begin{theorem}\label{at:G:key}
        Suppose $\Sigma$ is a $C^2$ embedded orientable hypersurface in $M$ with its unit normal vector $\vec{\nu}$. We denote by $\rH$ and $\rH_f$ the mean curvature of $\Sigma$ with respect to $\vec{\nu}$ and to $e^{-f}\vec{\nu}$ in $M$ and $M_f$, respectively. Then, for every $x$ in $\Sigma$, the following identity holds:
        $$
        \rH_{f}(x)=e^{-f}(\rH+n\langle \mathrm{D} f, \vec{\nu}\rangle )(x).
        $$
        Here $n$ is the dimension of $\Sigma$, $\langle\cdot,\cdot\rangle$ is the inner product of $M$, and $\mathrm{D} f$ is the gradient of $f$ in $M$.
    \end{theorem}

    \begin{proof}
        By definition, the mean curvature of $\Sigma$ is the divergence of $\vec{\nu}$ in $M$. Let $\vec{\nu}_f=e^{-f}\vec{\nu}$ be the unit normal vector to $\Sigma$ in $M_f$. Let $\dvol_f$ and $\dvol$ represent the volume form of $M_f$ and $M$, respectively. For any vector field $X$, the identity $\mathrm{d}(X\llcorner \dvol)=\div(X)\dvol$ is well-known, where $\div$ is the divergence of $M$.

        A direct computation yields:
        $$
        \begin{aligned}
            \rH_f e^{(n+1)f}\dvol &= \rH_f \dvol_f= \mathrm{d}(\vec{\nu}_f\llcorner \dvol_f)\\
            &=\mathrm{d}(e^{n f}\vec{\nu}\llcorner \dvol)\\
            &=e^{nf}(\rH+n\langle \mathrm{D} f, \vec{\nu}\rangle )\dvol,
        \end{aligned}
        $$
        where $\langle\cdot,\cdot\rangle$ is the inner product of $M$, and $\mathrm{D} f$ is the gradient of $f$ in $M$. Comparing both sides yields the conclusion.
    \end{proof}

    Recall that in our settings, $N$ is an $n$-dimensional $C^{2,\alpha}$ closed or compact Riemannian manifold with metric $\sigma$, and let $N^0$ be the interior of $N$. By a $C^{2,\alpha}$ compact Riemannian manifold, we mean a compact manifold with $C^{2,\alpha}$ boundary. A conformal product manifold $M_f$ is defined by
    \be \label{def:A:B}
    M_f:= \{N\PLH (a,b), e^{2f}(\sigma+\rd r^2)\}.
    \ene
    Here $f=f(x,r)$ is a smooth function in a neighborhood of $N\PLH(a,b)$.

    \begin{Rem}
        Recall that a warped product manifold is defined by
        $$
        \{   N\PLH (a,b), h^2(r)\sigma+\rd r^2\}
        $$
        where $h(r)\in C^\infty(a,b)$ is a positive function. By introducing a new parameter $s$ satisfying $\rd s=\frac{\rd r}{h(r)}$ and defining $f(s)=\ln h(r(s))$, we can express the metric of warped product manifolds in the form
        $$
        h^2(r)\sigma+\rd r^2=e^{2f(s)}(\sigma+\rd s^2).
        $$
        In view of the definition \eqref{def:A:B}, all warped product manifolds---including Euclidean spaces, hyperbolic spaces and hemispheres---are instances of conformal product manifolds.
    \end{Rem}

    By Theorem~\ref{at:G:key}, the mean curvature of the graph of a $C^2$ function $u$ with respect to the upward unit normal vector $e^{-f}(\F{-\nabla u+\partial_r}{\omega})$ in $M_f$ is
    \be \label{expression:LF}
    e^{-f}\Bigl(-\div\Bigl(\F{\nabla u}{\omega}\Bigr)+n\F{\langle \nabla u, \nabla f\rangle}{\omega}-n\F{\partial f}{\partial r}\F{1}{\omega} \Bigr),
    \ene
    where $\nabla u$ denotes the gradient of $u$ in $N$, $\nabla f$ is the gradient of $f$ for fixed $r$, $\div$ is the divergence of $N$, $\partial_r$ is the coordinate vector field on $(a,b)$, and $\omega=\sqrt{1+|\nabla u|^2}$.

    \begin{Rem}\label{nonatural:rk}
        By Theorem~\ref{at:G:key}, the graph of a $C^2$ function $u$ is minimal in $M_f$ if and only if $u$ solves the following PMC equation
        \be \label{eq:ms}
        -\div\Bigl(\F{\nabla u}{\omega}\Bigr)+n\F{\langle \nabla u, \nabla f\rangle}{\omega}-n\F{\partial f}{\partial r}\F{1}{\omega}=0.
        \ene
        This equation can also be expressed as
        \be \label{non:linear:PDE}
        \mathcal{L}_{\mathcal{F}}(u)=0, \quad \mathcal{L}_{\mathcal{F}}(u):=-\div\Bigl(\F{\nabla u}{\omega}\Bigr)-\mathcal{F}\Bigl(x,u, -\F{\nabla u}{\omega}, \F{1}{\omega}\Bigr)
        \ene
        with $\mathcal{F}(x,r, Y, t):=-n\langle Y, \nabla f\rangle +n\F{\partial f}{\partial r}t$ being a $C^{1,\alpha}$ function in the tangent bundle $T(N\PLH (a,b))$, where we adopt the notation $f=f(x,r)$. Given that the conformal factor $f$ can be chosen almost arbitrarily, the condition
        \be\label{condition:monotone}
        \F{\partial \mathcal{F}}{\partial r}\leq 0
        \ene
        is not naturally satisfied.
    \end{Rem}

    However, condition \eqref{condition:monotone} plays a crucial role in the interior gradient estimate and the Perron method for prescribed mean curvature equations, as evidenced in~\cite{GT01} and~\cite[Section~4]{zhou22a}. Indeed, this condition is essential for ensuring the existence of regular PMC graphs as follows.

    \begin{theorem}\label{thm:PDE:existence:}
        Suppose $N$ is a closed or $C^{2,\alpha}$ compact Riemannian manifold, and that $\mathcal{F}(x,r,Y,t):T(N\PLH \R)\rightarrow \R$ is a $C^{1,\alpha}$ function satisfying \eqref{condition:monotone}. Let $u_1,u_0$ be two functions in $C^2(N^0)\cap C(N)$ satisfying $u_1<u_0$ on $N^0$, and suppose that $(\rgr(u_1), \rgr(u_0))$ forms a barrier for $(N\PLH\R, \mathcal{F})$ in the sense of Definition~\ref{def:A}. Then there is a unique function $u$ in $C^{2,\alpha}(N^0)\cap C(N)$ that solves equation \eqref{non:linear:PDE} on $N^0$ and satisfies $u_1\leq u\leq u_0$.

        In the case where $N$ is closed, the boundary $\partial N$ is empty.
    \end{theorem}

    \begin{proof}
        By Definition~\ref{def:A}, we have $\mathcal{L}_{\mathcal{F}}(u_1) \leq 0$ and $\mathcal{L}_{\mathcal{F}}(u_0) \geq 0$. We apply the Perron method as in~\cite[Section~4]{zhou22a}. Under condition~\eqref{condition:monotone} and by~\cite[Theorem~A.2, Theorem~A.3]{zhou22a}, the argument in the proof of~\cite[Lemma~4.8]{zhou22a} implies that $u_1$ and $u_0$ are a Perron subsolution and a Perron supersolution of $\mathcal{L}_{\mathcal{F}}(u) = 0$, respectively.

        Now define the set
        $$
        \mathcal{S}_{u_0}=\{ u: \text{$u$ is a Perron subsolution of } \mathcal{L}_{\mathcal{F}}, u\in C(N), u \leq u_0\}.
        $$
        By~\cite[Lemma~4.13]{zhou22a} there exists a Perron solution $u \in C^{2,\alpha}(N^0)\cap C(N)$ satisfying $\mathcal{L}_{\mathcal{F}}(u)=0$ on $N^0$ and $u_1\leq u \leq u_0$. Uniqueness follows from~\eqref{condition:monotone} and the maximum principle.
    \end{proof}

    \begin{Rem}
        Theorem~\ref{thm:PDE:existence:} can be regarded as a generalization of earlier results in Euclidean spaces established by~\cite{BK74, TW83} and in hemispheres as in~\cite{HSW00}. In those works, the conditions imposed on PMC functions are special cases of condition~\eqref{condition:monotone}.
    \end{Rem}

    A barrier condition in the context of Riemannian manifolds is defined as follows.

    \begin{Def}\label{definiton:barrier:condition}
        Let $M^{n+1}$ be an $(n+1)$-dimensional compact manifold with smooth embedded boundary $\partial M$, let $\Gamma \subset \partial M$ be a closed smooth hypersurface (possibly empty), and let $\mathcal{H}$ be a $C^{1,\alpha}$ function on $TM$. Suppose $\partial M \backslash \Gamma=S_1\cup S_2$, where $S_1$ is connected, and both $S_1$ and $S_2$ are nonempty. We say that the triple $(S_1,S_2,\Gamma)$ is a barrier for $(M,\mathcal{H})$ if
        $$
        \rH_{S_1}(p)\geq \mathcal{H}(p,\vec{\nu}), \quad \rH_{S_2}(p)\leq  \mathcal{H}(p,\vec{\nu}),
        $$
        where $\rH_{S_i}:=\div(\vec{\nu})$ denotes the mean curvature of $S_i$ for $i=1,2$, and the unit normal vector $\vec{\nu}$ of $S_1$ ($S_2$) is outward (inward) with respect to $M$.
    \end{Def}

    For PMC hypersurfaces, the barrier condition described above plays a role analogous to boundary values of continuous functions on a closed interval. A PMC analogue of the intermediate value theorem can be stated as follows.

    \begin{theorem}
        Let $(S_1, S_2,\Gamma)$ be a barrier for $(M^{n+1},\mathcal{H})$ defined as above.
        \begin{enumerate}
            \item For $2\leq n\leq 7$, there exists a smooth embedded orientable hypersurface $\Sigma \subset M$ with boundary $\Gamma$ and mean curvature $\mathcal{H}$, as shown by Eichmair~\cite[Theorem~1.1]{Eich09}. Note that $\Gamma$ may be empty;
            \item For $n=2$, if $\Gamma$ is a null-homotopic rectifiable Jordan curve, $\mathrm{H}_2(M^3)\equiv 0$, and $\mathcal{H}\equiv \mu$ (a positive constant), then there is an embedded disk with constant mean curvature $\mu$ and boundary $\Gamma$, as established by Coskunuzer~\cite{Cos17}.
        \end{enumerate}
    \end{theorem}

    In Definition~\ref{definiton:barrier:condition}, $\partial M$ and $\Gamma$ are smooth and can be relaxed to piecewise smooth boundaries. In particular, the barrier condition in conformal product manifolds admits the following weaker form.

    \begin{Def}\label{def:A}
        Let $M_f$ be as in \eqref{def:A:B} and let $\mathcal{H}\in C^{1,\alpha}(T M_f)$. Given two functions $u_1,u_0 \in C^{2}(N^0)\cap C(N)$, we say that the pair $(\rgr(u_1),\rgr(u_0))$ is a barrier for $(M_f, \mathcal{H})$ if $u_1<  u_0$ on $N^0$, and the mean curvatures of $\rgr(u_1)$ and $\rgr(u_0)$ with respect to the upward-pointing unit normal vector satisfy
        $$
        \rH_{\rgr(u_1)}(p)\leq \mathcal{H}(p,\vec{\nu}),\quad
        \rH_{\rgr(u_0)}(p)\geq  \mathcal{H}(p,\vec{\nu}).
        $$
        Here $\rgr(u_i)$ is the graph of $u_i$, and $\rH_{\rgr(u_i)}$ is the mean curvature of $\rgr(u_i)$ for $i=1,0$. If $\partial N\neq \emptyset$, we further require that $u_1=u_0=\psi$ for some $\psi \in C^{1,\alpha}(\partial N)$.
    \end{Def}

    \subsection{Some facts on $\Lambda$-perimeter minimizers}
    In this subsection, we collect some properties on $\Lambda$-perimeter minimizers in an $(n+1)$-dimensional Riemannian manifold $M^{n+1}$. Let $\Omega \subset M$ be a fixed domain. Denote by $L^1(\Omega)$ the space of measurable functions $u$ on $\Omega$ such that $\int_{\Omega}u \dvol <\infty$, where $\dvol$ is the volume form of $M$. Let $\mathrm{div}$ and $\langle \cdot, \cdot \rangle$ denote the divergence and the Riemannian inner product on $M$, respectively. For any vector field $X$, let $\spt(X)$ denote the support of $X$, i.e., the closure of the set $\{x \in \Omega : X(x) \neq 0\}$. Denote by $\Gamma(TM)$ the set of all smooth vector fields in $TM$.

    \begin{Def}\label{def:BV}
        Let $u\in L^1(\Omega)$. We say that
        \begin{enumerate}
            \item $u$ has a weak derivative $g\in \Gamma(T M)$ if
            $$
            \int_{\Omega} u\, \div(X) \dvol=-\int_\Omega\langle X, g\rangle \dvol
            $$
            for any vector field $X\in \Gamma(TM)$ with $\spt(X)\subset\subset  \Omega$. The vector space of all functions $u$ with weak derivatives in $L^1(\Omega)$ is denoted by $W^{1,1}(\Omega)$. For a complete definition of $W^{k,p}(\Omega)$, see~\cite[Section~7.5]{GT01};
            \item $u$ is a function of bounded variation if the distributional derivative of $u$, $Du$, is a signed $TM$-valued Radon measure given by
            $$
            \int_{\Omega} u\, \div(X) \dvol=-\int_\Omega\langle X, \rd D u\rangle
            $$
            for any vector field $X\in \Gamma(TM)$, $\spt(X)\subset \subset \Omega$.
        \end{enumerate}
    \end{Def}

    \begin{Rem}
        It is well-known that $W^{1,1}(\Omega)\subset \BV(\Omega)$. When $u\in C^1(\Omega)$, the notations $\nabla u$ and $Du$ are used interchangeably to denote the gradient of $u$.
    \end{Rem}

    For any $TM$-valued measure $\nu$, we define its variation on $\Omega$ as
    $$
    |\nu|(\Omega)=\sup\Bigl\{\int_{\Omega}\langle X, \rd\nu\rangle : \langle X, X\rangle \leq 1, X\in \Gamma(TM), \spt(X)\subset \subset \Omega\Bigr\}.
    $$
    Then $u\in \BV(\Omega)$ if and only if $|D u|(\Omega)$ is finite. We say $u\in \BV_{loc}(\Omega)$ if $u\in \BV(W)$ for any bounded open set $W\subset \subset \Omega$.

    \begin{Def}[{\cite[P.\,122]{Mag12}}]\label{def:finite:locally:set}
        A measurable set $E$ is a set of locally finite perimeter or a Caccioppoli set if its characteristic function $\chi_E$ belongs to $\BV_{loc}(\Omega)$.

        The perimeter of $E$ in $\Omega$ is defined as
        $$
        P(E,\Omega)=|D\chi_E|(\Omega).
        $$
        The boundary of $E$, denoted by $\partial E$, is the set of all points $x\in \Omega$ such that
        $$
        0<\mathrm{vol}(B_r(x)\cap E)< \mathrm{vol}(B_r(x))
        $$
        for every $r>0$, where $B_r(x)$ denotes an embedded ball centered at $x$ with radius $r$.
    \end{Def}

    \begin{Rem}
        Notice that $P(E,\Omega)=\mathrm{H}^{n}(\partial^*E\cap \Omega)$, where $\mathrm{H}^{n}$ denotes the $n$-dimensional Hausdorff measure in $\Omega$ and $\partial^*E$ is the reduced boundary of $E$. In the case where $E$ has $C^1$ boundary, it holds that $P(E, \Omega)=\mathrm{H}^{n}(\Omega \cap \partial E)$.
    \end{Rem}

    In the following definition, we remove the restriction that $E\Delta F\subset \subset B_r(x_0)$ in the original definition in~\cite[Page~278]{Mag12}.

    \begin{Def}\label{Def:amb}
        Let $\Lambda\geq 0$. A Caccioppoli set $E$ is called a $\Lambda$-perimeter minimizer in an open set $W\subset M$ if for any Caccioppoli set $F$ satisfying $F\Delta E\subset  \subset W$, it holds that
        \be \label{eq:variation:hold}
        P(E, W)\leq P(F, W)+\Lambda \mathrm{vol}(E\Delta F).
        \ene
        Define the regular set $\mathrm{reg}(\partial E)$ as the set $\{x\in \partial E: \partial E \text{ is a } C^{1,\gamma}$ graph in a neighborhood of $x$ for some $\gamma\in (0,\F{1}{2})\}$ and the singular set $\mathrm{sing}(\partial E)$ as $\partial E\backslash \mathrm{reg}(\partial E)$.
    \end{Def}

    By~\cite[27.1]{Simon83}, for any Caccioppoli set $E$, let $[[\partial E]]$ denote the integer multiplicity current associated with $\partial E$, denoted by $\tau(\partial E, \xi)$. Here $\xi$ is the volume form of $\partial E$. The integer multiplicity of $[[\partial E]]$ is one almost everywhere. For any integer multiplicity current $T$ and an open set $W\subset M$, let $\mathbf{M}_W(T)$ denote the mass of $T$ in $W$. A useful fact is that for any Caccioppoli set $E$,
    $$
    \mathbf{M}_W([[\partial E]])= P(E, W).
    $$

    \begin{Def}[{\cite[Page~357]{DS93}}]
        We say $T$ is an $n$-dimensional $\Lambda$-minimizing integer multiplicity current in $W$ if
        $$
        \mathbf{M}_W(T)\leq \mathbf{M}_W(T+\partial Q)+\Lambda \mathbf{M}_W(Q)
        $$
        holds for any $(n+1)$-dimensional integer multiplicity current $Q$ with $\spt(Q)\subset \subset W$.
    \end{Def}

    An example of a $\Lambda$-perimeter minimizer is constructed as follows.

    \begin{theorem}\label{thm:lambda:perimeter}
        Let $\Omega$ be a bounded $C^2$ domain in $M^{n+1}$. Suppose $u$ is a bounded function in $C^2(\Omega)$ with Lipschitz boundary value $\psi$ on $\partial \Omega$ satisfying
        \be\label{condition:mean:curvature}
        \sup_{x\in \Omega } \Bigl\{\Bigl|\div\Bigl(\F{\nabla u}{\sqrt{1+|\nabla u|^2}}\Bigr)\Bigr|\Bigr\}\leq \kappa
        \ene
        for some constant $\kappa$. Let $T$ be the integer multiplicity current $[[\{(x, u(x)):x\in \Omega \}]]$ and define $\Omega_a =\{x\in M: d(x,\Omega)<a\}$ for any positive $a$. Then there exist positive constants $\Lambda:=\Lambda(\Omega,\kappa)$ and $\delta:=\delta(\Omega, \kappa)$ such that
        \begin{enumerate}
            \item The subgraph $\mathrm{sub}(u)=\{(x,t):x\in \partial\Omega, t<u(x)\}$ is a $\Lambda$-perimeter minimizer in $\Omega_\delta\PLH \R\backslash \rgr(\psi)$ with respect to the product metric;
            \item $T$ is an $(n+1)$-dimensional $\Lambda$-minimizing integer multiplicity current in $\Omega_\delta \PLH \R$ with boundary $\partial T=\rgr(\psi)$.
        \end{enumerate}
    \end{theorem}

    \begin{Rem}\label{Density:T}
        Suppose $\psi$ belongs to $C^{1,\alpha}(\partial\Omega)$. Observe that for any $a\in \rgr(\psi)$, the density satisfies $\Theta(\|T\|, a)=\F{1}{2}$ because in a sufficiently small neighborhood of $\rgr(\psi)$, the support $\spt(T)$, which coincides with $\rgr(u)$, forms only one connected component. By the result of Duzaar and Steffan~\cite[Theorem~4.5\,(1)]{DS93}, in a neighborhood of $\rgr(\psi)$, the support $\spt(T)$ is a $C^{1,\gamma'}$ embedded manifold-with-boundary for any $\gamma'\in (0, \F{\alpha}{2})$.
    \end{Rem}

    \begin{proof}
        Item (1) has been shown in~\cite[Lemma~2.12]{zhou22a}. Since $\rgr(\psi)$ is of class $C^{1,\alpha}$, it follows directly that $\partial T=\rgr(\psi)$. By definition, the mass of $T$ equals $P(\mathrm{sub}(u),\Omega\PLH\R)$.

        Now we proceed to the proof of item (2).
        Let $V$ be a bounded domain with $C^2$ boundary and let $\delta>0$ be sufficiently small. Define
        $V_\delta=\{x \in M : d(x,V) < \delta\}$, $(\partial V)_\delta=\{x\in M : d(x, \partial V) < \delta\}$, and $V^\delta=\{x\in V : d(x, \partial V) > \delta\}$.

        For any point $p=(x,t) \in M \PLH \R$, define the signed distance function:
        $$d(p)=\mathrm{sign}(x) d(x, \partial \Omega)$$
        where $\mathrm{sign}(x)=1$ for $x \in \Omega$ and $\mathrm{sign}(x)=-1$ for $x \in M \backslash \Omega$.

        For sufficiently small $\delta$, the function $d$ is $C^{2}$ in a neighborhood of $\partial\Omega\PLH\R$. The $(n+1)$-form $\sigma_0:=D d\llcorner \dvol$ is then well-defined near $\partial\Omega\PLH\R$. Let $\sigma_1:=\vec{\nu}\llcorner \dvol$, where $\vec{\nu}$ is the vertical translation of the normal vector of the graph of $u$ in $\Omega\PLH\R$.

        Differentiating $\sigma_0$ and $\sigma_1$, we find
        $$\rd\sigma_0=f_0 \dvol \quad \text{ near}\quad  (\Omega_\delta\backslash \Omega)\PLH\R,\quad \rd\sigma_1= f_1 \dvol\quad \text{in}\quad\Omega\PLH\R, $$
        with $|f_i| \leq \Lambda := \Lambda(\Omega, \kappa)$ for $i = 0, 1$, provided $\delta = \delta(\Omega, \kappa)$ is chosen sufficiently small, as implied by~\eqref{condition:mean:curvature}.

        Let $W$ be an open bounded set in $\Omega_\delta\PLH\R$. For any $\varepsilon>0$, let $X_\varepsilon$ be any smooth $(n+1)$-form on $M\PLH \R$ satisfying
        \begin{enumerate}
            \item $\langle X_\varepsilon, X_\varepsilon \rangle \le 1$, with $W^\varepsilon\subset \spt(X_\varepsilon) \subset \subset W$;
            \item $X_\varepsilon=\sigma_1$ in $W^\varepsilon \cap (\Omega^\varepsilon\PLH \R)$;
            \item $X_\varepsilon=\sigma_0$ in $W^\varepsilon \cap ((\partial \Omega)_\varepsilon\PLH\R)$.
        \end{enumerate}
        Fix any Caccioppoli set $E\subset \subset W$ and any multiplicity function $\theta$. Let $Q$ be the $(n+2)$-dimensional integer multiplicity current $\tau(E,\theta)$. Define:
        $$\theta_0=\chi_{(\Omega_\delta\backslash \Omega)\PLH\R}\theta,\quad \theta_1=\chi_{\Omega\PLH\R}\theta,$$
        where $\chi_E$ denotes the characteristic function of $E$. Then $Q=Q_0+Q_1$, where $Q_0=\tau (E, \theta_0)$ and $Q_1=\tau(E,\theta_1)$.

        Let $T^a$ be the integer multiplicity current $\partial [[\mathrm{sub}(u)]]$. Then $T^a=T+T_0$, where $T_0$ is the integer multiplicity current $[[\mathrm{sub}(u|_{\partial\Omega})]]$ in $\partial\Omega\PLH\R$. It follows that
        $$T^a(X_\varepsilon)=(T+\partial Q)(X_\varepsilon)+T_0(X_\varepsilon)-\partial Q(X_\varepsilon).
        $$
        By the definitions of $\sigma_0$, $\sigma_1$, and the smoothness of $X_\varepsilon$, $\spt(T)$, and $\spt(T_0)$, we have:
        $$
        \lim_{\varepsilon\rightarrow 0} T^a(X_\varepsilon)=\mathbf{M}_W(T)+\mathbf{M}_W(T_0).
        $$

        Since $\mathbf{M}_W(Q_i) = \int_E |\theta_i|$ for $i = 0, 1$, we have $\mathbf{M}_W(Q_0) + \mathbf{M}_W(Q_1) = \mathbf{M}_W(Q)$. Moreover, by~\cite[(2.19)]{zhou22a} and the definitions of $\sigma_0$ and $\sigma_1$,
        $$
        \lim_{\varepsilon\rightarrow 0} \partial Q(X_\varepsilon)=\lim_{\varepsilon\rightarrow 0} (Q_0+ Q_1)(\rd X_\varepsilon)=Q_0(\rd \sigma_0)+Q_1(\rd\sigma_1)\geq -\Lambda \mathbf{M}_{W}(Q).
        $$
        Thus,
        $$
        \begin{aligned}
            \mathbf{M}_W(T)+\mathbf{M}_W(T_0) &= \lim_{\varepsilon\rightarrow 0}T^a(X_\varepsilon)\\
            &\leq  \mathbf{M}_W(T+\partial Q)+\mathbf{M}_W(T_0)+\Lambda \mathbf{M}_W(Q),
        \end{aligned}
        $$
        which completes the proof.
    \end{proof}

    We record the following regularity result.

    \begin{theorem}[{\cite[Theorem~28.1]{Mag12}}]\label{thm:regularity}
        Let $M$ be an $(n+1)$-dimensional Riemannian manifold. Suppose $E$ is a $\Lambda$-perimeter minimizer in an open set $W\subset M$. Then
        \begin{enumerate}
            \item if $2\leq n\leq 6$, $\mathrm{sing}(\partial E)$ is empty;
            \item if $n=7$, $\mathrm{sing}(\partial E)$ is a discrete set;
            \item if $n\geq 8$, $\mathrm{H}^s(\mathrm{sing}(\partial E))=0$ for any $s>n-7$.
        \end{enumerate}
    \end{theorem}

    We say a sequence of Caccioppoli sets $\{E_i\}_{i=1}^\infty$ locally converges to a Caccioppoli set $E$ in an open set $\Omega\subset M$ if for any Borel set $A\subset \subset \Omega$
    $$
    \lim_{i\rightarrow +\infty}	 \int_A |\chi_E-\chi_{E_i}|\dvol=0,
    $$
    where $\dvol$ is the volume form of $M$. Here $\chi_E$ is the characteristic function of $E$.

    The following properties of $\Lambda$-perimeter minimizers are very useful.

    \begin{lemma}\label{lm:slice:convergence}
        Let $\{E_j\}_{j=1}^\infty$ be a sequence of $C^2$ $\Lambda$-perimeter minimizers in a bounded open set $W\subset M$. Then:
        \begin{enumerate}
            \item ({\cite[Proposition~21.13  Theorem~21.14]{Mag12}}) There exists a $\Lambda$-perimeter minimizer $E$ such that a subsequence of $\{E_j\}_{j=1}^\infty$, still denoted by $\{E_j\}_{j=1}^\infty$, converges locally to $E$ in $W$;
            \item Suppose in addition that $\partial E$ is regular in $W$, i.e., $\mathrm{sing}(\partial E) \cap W = \emptyset$. Then:
            \begin{enumerate}
                \item the unit normal vectors of $\{\partial E_j\}_{j=1}^\infty$ and $\partial E$ are locally uniformly continuous with respect to the distance in $W$;
                \item $\{\partial E_j\}_{j=1}^\infty$ converges locally uniformly to $\partial E$ in the $C^{1,\gamma}$ sense in $W$ for any $\gamma \in (0,\F{1}{2})$.
            \end{enumerate}
        \end{enumerate}
    \end{lemma}

    \begin{proof}
        We provide only a sketch of the proof for item (2). This result is an enhanced version of~\cite[Theorem~26.6]{Mag12}, which is implicitly contained in~\cite[Lemma~3.6]{zhou22a} together with the Allard regularity theorem. Without loss of generality, we assume that $W$ is isometrically embedded into some Euclidean space $\mathbb{R}^{n+k}$.

        Let $K$ be a compact subset of $W$.
        By the regularity of $\partial E$,~\cite[Lemma~3.6]{zhou22a} and the Allard regularity theorem (see, e.g.,~\cite[Theorem~3.5]{zhou22a},~\cite{Simon83}), there is a sufficiently small $\rho>0$ such that for any $E_i$ in $\{E_i\}_{i=1}^\infty$, any $x\in  K\cap \partial E_i$, any $\gamma \in (0,\F{1}{2})$, one can find constants
        $\delta=\delta(n,k,p)$, $\gamma_0= \gamma_0(n,k,\gamma)\in (0,1)$, a linear isometry $q\in \R^{n+k}$ and $u\in C^{1,\gamma}(\bar{B}^n_{\gamma_0\rho}(0),\R^{k-1})$ satisfying
        $$u(0)=0, \quad x=q((0,u(0))),\quad \partial E_i\cap B^{n+k}_{\gamma_0\rho}(x)=q(\rgr(u))\cap B^{n+k}_{\gamma_0\rho}(x)$$
        and
        \be \label{equality}
        \rho^{-1}\sup_{B^n_{\gamma_0\rho}(0)} |u|+\sup_{B^n_{\gamma_0\rho}(0)}|\nabla u|+\rho^{\gamma }\sup_{\substack{y',y\in B^n_{\gamma_0\rho}(0),\\
                y'\neq y}}\F{|\nabla u(y')-\nabla u(y)|}{|x-y|^{\gamma}}\leq c\delta^{\F{n}{4}}
        \ene
        where $c=c(n,k,\gamma)$.
        This implies that for every $E_i \in \{E_j\}_{j=1}^\infty$, the boundary $\partial E_i \cap K$ is a locally uniformly $C^{1,\gamma}$ graph. Consequently, their unit normal vectors are locally uniformly continuous with respect to the distance in $W$. By the compactness property of $\Lambda$-minimizers, the same conclusion holds for $E$. This establishes item (a).

        Furthermore, by~\cite[Theorem~26.6]{Mag12}, if a sequence $\{x_j\}_{j=1}^\infty$ with $x_j \in \partial E_j$ converges to a point $x \in \mathrm{reg}(\partial E)$, then the corresponding unit normal vectors also converge. Combining this with the result of item (a), we obtain item (b).
    \end{proof}

    \section{Limit of bounded graphs with bounded mean curvature}
    In this section we consider the $L^1$ limit of a sequence of $C^2$ functions whose graphs have bounded mean curvature.

    Throughout this section, $N$ denotes a fixed Riemannian manifold. For any $x\in N$, let $B_r(x)$ denote the embedded ball centered at $x$ with radius $r$ in $N$ and let $\mathrm{vol}(B_r(x))$ denote its volume. We now adapt the definitions from~\cite[Definition~3.70]{AFP00} to Riemannian manifolds.

    \begin{Def}\label{def:ap}
        Let $\Omega$ be an open subset of $N$, and let $u \in L^1(\Omega)$. A number $z\in \R$ is said to be an approximate limit of $u$ at $x\in \Omega$ if
        \be \label{eq:property}
        \lim_{r\rightarrow 0}\F{1}{\mathrm{vol}(B_r(x))}\int_{B_r(x)}|u(y)-z|\dvol=0.
        \ene
        The set $S_u \subset \Omega$ where this property fails is called the \emph{approximate discontinuity set}. For any $x\in \Omega\backslash S_u$, the uniquely determined value $z$ in~\eqref{eq:property} is called the approximate limit of $u$ at $x$ and is denoted by $\tilde{u}(x)$.

        If $N\equiv \R^n$, for any $x\in \Omega\backslash S_u$, we say $u$ is approximately differentiable at $x$ if there is a vector $L\in T_x \R^n$ such that
        \be \label{eq:property:B}
        \lim_{r\rightarrow 0}\F{1}{\mathrm{vol}(B_r(x))}\int_{B_r(x)}\frac{|u(y)-\tilde{u}(x)-\langle L, (y-x)\rangle|}{r}\dvol=0.
        \ene
        Then such $L$, uniquely determined by~\eqref{eq:property:B}, is denoted by $\nabla u(x)$.

        When $N$ is a Riemannian manifold, for any $x\in \Omega$, we say $u$ is approximately differentiable at $x$ if for any diffeomorphism $\phi^{-1}:\Omega \rightarrow \R^n$, $u\circ \phi$ is approximately differentiable at $\phi^{-1}(x)$. Moreover, there is a unique tangent vector $L$ in $T_x N$ such that
        $$
        \nabla (u\circ \phi )(\phi^{-1}(x))=\phi_* (L)(\phi^{-1}(x))
        $$
        where $\phi_*:T_x N \rightarrow T_{\phi(x)}\R^n$ is the derivative of $\phi$. In this case $L$ is written as $\nabla u(x)$.
    \end{Def}

    The following result of Calder\'on-Zygmund type holds.

    \begin{theorem}[{\cite[Theorem~3.83]{AFP00}}]\label{Thm:ap}
        Suppose $u\in \BV(\Omega)$. Then $u$ is approximately differentiable $\mathcal{L}^n$-a.e. Moreover, $\nabla u$ is the density vector of $D u$ with respect to the Lebesgue measure $\dvol$.
    \end{theorem}

    By the Lebesgue decomposition Theorem, for any BV function $u\in \BV(\Omega)$, it holds that
    \be \label{eq:Lesbegue:decomposition}
    D u=\nabla u \dvol+ D^s u
    \ene
    where $D^su$ is the singular part of $D u$ with respect to the Lebesgue measure $\dvol$.

    \begin{lemma}\label{lemma:A}
        Suppose $u\in \BV(\Omega)$. Then $u\in W^{1,1}(\Omega)$ if and only if $|D^su|(\Omega)\equiv 0$.
    \end{lemma}

    For any $u\in \BV(\Omega)$, define a functional
    \be\label{def:area:functional}
    \begin{split}
        \int_{\Omega}\sqrt{1+|D u|^2}:&=\sup\Bigl\{\int_{\Omega}h+u \div(X)\dvol: \\
        & \spt(h)\cup \spt(X)\subset \subset \Omega, h^2+\langle X, X\rangle\leq 1\Bigr\}
    \end{split}
    \ene
    where the supremum is taken over all smooth functions $h$ and smooth vector fields $X$ on $\Omega$.

    We continue to use $P$ to denote the perimeter of a Caccioppoli set in the product manifold $N\PLH \R$ when no confusion arises.
    The following result is well-known.

    \begin{theorem}[{\cite[Theorem~14.6]{Giu84}}]\label{thm:AD}
        Suppose $u\in \BV(\Omega)$ and $\Omega$ is open and bounded. Then
        $$
        \begin{aligned}
            P(\mathrm{sub}(u),\Omega\PLH \R) &=\int_{\Omega}\sqrt{1+|D u|^2}\\
            &=\int_{\Omega} \sqrt{1+|\nabla u|^2}\dvol+|D^su|(\Omega)
        \end{aligned}
        $$
        where $\nabla u(x)$ is defined by~\eqref{def:ap}.
    \end{theorem}

    \begin{proof}
        The first equality is given in~\cite[Theorem~14.6]{Giu84}. We now prove the second equality. By~\eqref{eq:Lesbegue:decomposition} and item (2) in Definition~\ref{def:BV}, we have
        $$
        \int_{\Omega}h+u \div(X)\dvol=\int_{\Omega}(h-\langle X, \nabla u\rangle) \dvol-\int_{\Omega}\langle X, \mathrm{d} D^su\rangle.
        $$
        The right-hand side is bounded above by $\int_{\Omega} \sqrt{1+|\nabla u|^2}\dvol+|D^su|(\Omega)$. The equality follows from the fact that the support of the measure $|D^su|$ is a set of Lebesgue measure zero.
    \end{proof}

    Combining~\eqref{def:area:functional} with the above theorem, we obtain the following result.

    \begin{cor}\label{local:bound}
        Let $u\in \BV(\Omega)$. Then
        $$
        P(\mathrm{sub}(u),\Omega\PLH\R)\geq \max\{\mathrm{vol}(\Omega), |D u|(\Omega)\}.
        $$
    \end{cor}

    The compactness of Radon measures gives the following compactness result for BV functions.

    \begin{theorem}\label{thm:compact}
        Suppose $\{u_k\}_{k=1}^\infty$ is a sequence in $\BV(\Omega)$ satisfying
        $$
        \sup_{k \ge 1}\Bigl\{\int_{\Omega }|u_k|\dvol+|D u_k|(\Omega)\Bigr\}\leq \kappa,
        $$
        for some positive $\kappa$. Then there exists a subsequence of $\{ u_k\}_{k=1}^\infty$ converging to $u$ in $ L_{loc}^1(\Omega)$.
    \end{theorem}

    The following property on $\Lambda$-perimeter minimizers is particularly useful.

    \begin{theorem}\label{thm:open:regularity}
        Let $v\in \BV(\Omega)$ be a bounded function such that $\mathrm{sub}(v)$ is a $\Lambda$-perimeter minimizer in $\Omega\PLH\R$. Then $v$ is of class $C^{1,\gamma}$ for any $\gamma \in (0,\F{1}{2})$ over a dense open set $\Omega_{v}\subset \Omega$. In addition, $\mathrm{vol}(\Omega_v)=\mathrm{vol}(\Omega)$.
    \end{theorem}

    \begin{proof}
        By~\cite[Theorem~26.3, 28.1]{Mag12}, there exists a closed set $S\subset \partial \mathrm{sub}(v)$ such that $ \mathrm{H}^s(S)=0$ for any $s>n-7$, and $\partial \mathrm{sub}(v)\backslash S$ is $C^{1,\gamma}$ orientable and embedded for some $\gamma\in (0,\F{1}{2})$. Let $S'$ be the vertical projection of $S$ from $\Omega\PLH\R$ into $\Omega$, and define $\Omega_0=\Omega \backslash S'$.
        Then
        $$\mathrm{vol}(\Omega\backslash \Omega_0)\leq \mathrm{H}^{n-1}(S)=0,$$
        and for each $x \in \Omega_0$ and $p=(x,v(x))$, there exists a neighborhood $U\subset M \PLH\R$ of $p$ such that the boundary of $\mathrm{sub}(u)$ is $C^{1,\gamma}$ in $U$.

        Now define the set
        \be \label{def:omega1}
        \Omega_v=\{ x\in \Omega_0: \nabla v\text{ exists and is finite}\}
        \ene
        where $\nabla v$ is the approximate differential of $v$, as in~\eqref{eq:Lesbegue:decomposition} with the decomposition
        $$
        D v=\nabla v \dvol +D^s v.
        $$

        By Theorem~\ref{Thm:ap}, it follows that
        \be \label{a:characteristic}
        \mathrm{vol}(\Omega\backslash \Omega_v)=0.
        \ene
        Fix any $x\in \Omega_v$ and let $p=(x,v(x))$. Since $\mathrm{sub}(v)$ is a $\Lambda$-perimeter minimizer, $\partial \mathrm{sub}(v)$ is embedded, connected, and $C^{1,\gamma}$ in a neighborhood of $p$. Its unit normal vector at $p$ in $\Omega\PLH\R$ is given by
        $$\vec{\nu}(p) =\F{\partial_r-\nabla v}{\sqrt{1+|\nabla v|^2}},$$
        where $\partial_r$ is the tangent vector of $\R$ pointing to positive infinity. Moreover, $\langle \vec{\nu}_p,\partial_r\rangle>0$.

        Since the unit normal vector $\vec{\nu}$ is continuous near $p$, it follows that $\langle \vec{\nu},\partial_r\rangle >0$ in a sufficiently small neighborhood of $p$. As a result, near $p$, the boundary $\mathrm{sub}(v)$ is a $C^{1,\gamma}$ graph over an open set containing $x$ in $\Omega$.

        Since $x\in \Omega_v$ is arbitrary, we conclude that $\Omega_v$ is open in $\Omega$. From~\eqref{a:characteristic}, $\Omega_v$ is dense in $\Omega$, which completes the proof.
    \end{proof}

    As an application of Theorem~\ref{thm:open:regularity}, the following result provides motivation for Theorem~\ref{key:thm}.

    \begin{lemma}\label{lm:fullmeasure}
        Let $v\in W^{1,1}(\Omega)$ and let $\mathrm{sub}(v)$ be a $\Lambda$-perimeter minimizer in $\Omega\PLH \R$. Let $\Sigma =\partial \mathrm{sub}(v)\cap (\Omega\PLH \R)$. Then there is an open dense set $\Omega_v$ such that $\Sigma \cap (\Omega_v\PLH\R)$ is $C^{1,\gamma}$ and dense in $\Sigma$ with respect to the $n$-dimensional Hausdorff measure, where $\dim\Omega=n$.
    \end{lemma}

    \begin{proof}
        Let $\Omega_v\subset \Omega$ be the open dense set obtained as in Theorem~\ref{thm:open:regularity}. Let $\mathrm{H}^n$ denote the $n$-dimensional Hausdorff measure in $\Omega\PLH\R$ with respect to the product metric.

        Observe that $\Sigma$ in $\Omega_v\PLH \R$ is a $C^{1,\gamma}$ graph. Then
        \begin{equation}\label{eq:vpart_zero}
            \begin{split}
                \mathrm{H}^n((\Omega_v\PLH \R)\cap \Sigma)&=\int_{\Omega_v}\sqrt{1+|\nabla v|^2}\dvol\\
                &=\int_{\Omega}\sqrt{1+|\nabla v|^2}\dvol= \mathrm{H}^n(\Sigma).
            \end{split}
        \end{equation}
        The second equality follows from $v\in W^{1,1}(\Omega)$ and $\mathrm{H}^n(\Omega\backslash \Omega_v)=0$.

        Therefore, we obtain the conclusion.
    \end{proof}

    \begin{theorem}\label{key:thm}
        Suppose $\Omega$ is a bounded $C^2$ domain.
        Define $\mathcal{D}_{\kappa}(\Omega)$ as the $L^1(\Omega)$ closure of $C^2$ functions $u$ on $\Omega$ satisfying
        \be\label{condition:A}
        \max_{\bar{\Omega}}\{|u|, |\div(\F{Du}{\sqrt{1+|Du|^2}})|\}\leq \kappa
        \ene
        for some positive constant $\kappa>0$. Then $\mathcal{D}_{\kappa}(\Omega)\subset W^{1,1}(\Omega)$ and is a compact set in $L^1(\Omega)$. For each $v\in \mathcal{D}_{\kappa}(\Omega)$, there exists an open dense set $\Omega_v$ such that $v$ is of class $C^{1,\gamma}$ on $\Omega_v$ for any $\gamma\in (0, \F{1}{2})$.
    \end{theorem}

    \begin{proof}
        Let $\{u_k\}_{k=1}^\infty$ be an arbitrary $C^2$ sequence satisfying~\eqref{condition:A}. By Theorem~\ref{thm:lambda:perimeter}, each $\mathrm{sub}(u_k)$ is a $\Lambda$-perimeter minimizer for some fixed $\Lambda$ and $|u_k|\leq \kappa$. Then
        $$
        \begin{aligned}
            P(\mathrm{sub}(u_k),\Omega\PLH \R)&\leq P(\mathrm{sub}(u_k)\cup \{\Omega \PLH (-2\kappa, 2\kappa)\},\Omega\PLH \R)\\
            &\quad+4\Lambda \kappa\mathrm{vol}(\Omega).
        \end{aligned}
        $$
        Note that
        $$
        P(\mathrm{sub}(u_k)\cup \{\Omega\PLH (-2\kappa, 2\kappa)\},\Omega\PLH \R)\leq 4\kappa \mathrm{H}^{n-1}(\partial \Omega)+2\mathrm{vol}(\Omega),
        $$
        where $n=\dim N$. As a result,
        \be\label{def:perimeter:bound}
        P(\mathrm{sub}(u_k), \Omega\PLH (-2\kappa, 2\kappa))\leq C
        \ene
        for some positive constant $C$ independent of $k$.
        By Corollary~\ref{local:bound}, it follows that
        \be \label{def:bound}
        \int_{\Omega}|u_k|\dvol + |D u_k|(\Omega)\leq C+\kappa \mathrm{vol}(\Omega).
        \ene
        By the compactness of BV functions, there exists a function $v\in \BV(\Omega)$ such that $u_k\rightarrow v$ in $L^1(\Omega)$. Therefore, $\mathcal{D}_\kappa(\Omega)$ is compact in $L^1(\Omega)$. By~\eqref{condition:A} and Lemma~\ref{lm:slice:convergence}, it follows that $\mathrm{sub}(v)$ is a $\Lambda$-perimeter minimizer in $\Omega\PLH\R$. Therefore, by Theorem~\ref{thm:open:regularity}, $v$ is of class $C^{1,\gamma}$ on an open and dense set. It remains to show that $\mathcal{D}_{\kappa}(\Omega)\subset W^{1,1}(\Omega)$.

        Since $\mathrm{sub}(v)$ is a $\Lambda$-perimeter minimizer in $\Omega\PLH\R$, Theorem~\ref{thm:regularity} implies the existence of a closed set $S\subset \partial \mathrm{sub}(v)$ such that $\mathrm{H}^{n-2}(S)=0$ and $\partial \mathrm{sub}(v)\backslash S$ is a $C^{1,\gamma}$, embedded and orientable hypersurface, where $n=\dim\Omega$. Define $\Omega_S=\Omega\backslash \pi (S)$ with $\pi :\Omega \PLH \R\rightarrow \Omega$ given by $\pi(x,r)=x$. Then $\Omega_S$ is open and dense in $\Omega$. Since $\mathrm{H}^{n-1}(\pi(S)\PLH\R)=0$, we have
        \be \label{fact:key:center}
        P(\mathrm{sub}(v),\Omega_S\PLH \R)=P(\mathrm{sub}(v),\Omega\PLH\R), \text{ and } |D^sv|(\pi(S))=0,
        \ene
        where $D^s(v)$ is defined in~\eqref{eq:Lesbegue:decomposition}.

        Now fix any open set $W\subset \subset \Omega_S$. By conclusion (2) of Lemma~\ref{lm:slice:convergence}, $\partial \mathrm{sub}(u_k)$ converges uniformly to $\partial \mathrm{sub}(v)$ in $W\PLH (-\kappa, \kappa)$ in the $C^{1,\gamma}$ sense. This implies:
        \begin{enumerate}
            \item[(i)] We have
            $$ \label{fact:subu_k}
            \lim_{k\rightarrow \infty }P(\mathrm{sub}(u_k),W\PLH\R)=P(\mathrm{sub}(v),W\PLH\R).
            $$
            \item [(ii)] $\{\nabla u_k\}_{k=1}^\infty$ converges to $\nabla v$ uniformly on any closed subset of $W\cap \Omega_v$, where $\Omega_v$ is defined as in~\eqref{def:omega1}.
        \end{enumerate}
        By Theorem~\ref{thm:open:regularity}, $v$ is $C^{1,\gamma}$ on the open dense set $\Omega_v$. Combining Theorem~\ref{thm:AD} and~\eqref{fact:subu_k} together, one sees that
        \be \label{ND}
        \begin{split}
            \lim_{k\rightarrow \infty}\int_{W}\sqrt{1+|\nabla u_k|^2}\dvol =\int_{W}\sqrt{1+|\nabla v|^2} \dvol+|D^sv|(W).
        \end{split}
        \ene

        Now we show that $|D^sv|(W)=0$. To this end, fix any $\varepsilon>0$ and let $\Omega_v$ be as in~\eqref{def:omega1}. Define $F_{v,s}$ as the closed set $\{x\in W\cap \Omega_v: |\nabla v|\leq s\}$ such that
        \be \label{mid:step}
        P(\mathrm{sub}(v), (W\cap \Omega_v\backslash F_{v,s})\PLH\R)=\int_{W\cap \Omega_v\backslash F_{v,s}}\sqrt{1+|\nabla v|^2}\dvol <\F{\varepsilon}{8}.
        \ene
        Since $W\cap \Omega_v\backslash F_{v,s}$ is open in $\Omega_S$, by (i) and~\eqref{mid:step}, there exists a $k_0>0$ such that for any $k>k_0$,
        \be \label{mid:step:one}
        \int_{(W\cap \Omega_v\backslash F_{v,s})\PLH\R}\sqrt{1+|\nabla u_k|}\dvol=P(\mathrm{sub}(u_k), (W\cap \Omega_v\backslash F_{v,s})\PLH\R)<\F{\varepsilon}{4}.
        \ene
        By fact (ii), there exists a $k_1$ such that for all $k>k_1$
        \be \label{mid:step:two}
        \Bigl|\int_{F_{v,s}}\sqrt{1+|\nabla u_k|}\dvol -\int_{F_{v,s}}\sqrt{1+|\nabla v|^2}\dvol\Bigr|\leq \F{\varepsilon}{2}.
        \ene
        For all $k\geq \max\{k_0,k_1\}$, combining~\eqref{mid:step},~\eqref{mid:step:one} and~\eqref{mid:step:two} yields
        \be
        \Bigl|\int_{W\cap \Omega_v}\sqrt{1+|\nabla u_k|^2}\dvol-\int_{W\cap \Omega_v}\sqrt{1+|\nabla v|^2}\dvol\Bigr|\leq \varepsilon.
        \ene
        Hence,
        $$
        \lim_{k\rightarrow \infty}\int_{W\cap \Omega_v}\sqrt{1+|\nabla u_k|^2}\dvol=\int_{W\cap \Omega_v}\sqrt{1+|\nabla v|^2}\dvol.
        $$
        Since $\Omega_v$ is open dense in $\Omega$ and each $u_k$ is $C^2$, we get
        $$
        \lim_{k\rightarrow \infty}\int_{W}\sqrt{1+|\nabla u_k|^2}\dvol=\int_{W \cap \Omega_v}\sqrt{1+|\nabla v|^2}\dvol.
        $$
        From~\eqref{ND}, it follows that $|D^sv|(W)=0$. Since $W$ is arbitrary, $|D^sv|(\Omega_S)=0$.

        From~\eqref{fact:key:center}, we have $|D^sv|(\Omega)=0$. By Lemma~\ref{lemma:A}, $v\in W^{1,1}(\Omega)$, which completes the proof.
    \end{proof}

    We conclude this section by establishing $C^{1,\gamma}$ convergence for $C^2$ bounded graphs with bounded prescribed mean curvature and fixed boundary. Here, a compact $C^{1,\alpha}$ manifold refers to one with boundary of class $C^{1,\alpha}$.

    \begin{theorem}\label{thm:boundary:convergence}
        Let $\Omega\subset N$ be a bounded $C^2$ domain with $\dim N=n$. Suppose $\{u_k\}_{k=1}^\infty$ is a sequence in $C^2(\Omega)\cap C(\bar{\Omega})$ satisfying~\eqref{condition:A} for some $\kappa>0$, such that $u_k=\psi$ on $\partial\Omega$ for some $\psi$ in $C^{1,\alpha}(\partial\Omega)$ and some $\alpha\in (0,1)$. Then, for any number $\gamma\in (0,\F{\alpha}{2})$, the following hold:
        \begin{enumerate}
            \item In a neighborhood of $\mathrm{gr}(\psi)$, the sequence $\{\mathrm{gr}(u_k)\}_{k=1}^\infty$ forms a family of uniformly compact $C^{1,\gamma}$ manifolds with boundary $\mathrm{gr}(\psi)$;
            \item If in addition $2\leq n\leq 7$, then a subsequence of $\{\mathrm{gr}(u_k)\}_{k=1}^\infty$ converges uniformly to a compact $C^{1,\gamma}$ embedded hypersurface $\Sigma$ with boundary $\mathrm{gr}(\psi)$. As a result, $\Sigma$ is homeomorphic to the closure of $\Omega$.
        \end{enumerate}
    \end{theorem}

    \begin{proof}
        Without loss of generality, we assume that $N\PLH \R$ is isometrically embedded into a Euclidean space $\R^{n_0}$.

        Let $T_k$ denote the integer multiplicity current induced by $\mathrm{gr}(u_k)$, let $\mathrm{gr}(\psi)$ denote the integer multiplicity current associated with the graph of $\psi$. By the compactness of integer multiplicity currents and Theorem~\ref{thm:lambda:perimeter}, a subsequence of $\Lambda$-minimizing integer multiplicity currents $\{T_k\}_{k=1}^\infty$ converges to a $\Lambda$-minimizing integer multiplicity current $T$, with $\partial T=\mathrm{gr}(\psi)$. Let $\Sigma$ be the support of $T$.

        As in~\eqref{def:bound}, there exists a $v\in W^{1,1}(\Omega)$ such that $\{u_k\}_{k=1}^\infty$ converges to $v$ in $L^1(\Omega)$. By Lemma~\ref{lm:slice:convergence}, $\Sigma$ is the set $\partial \mathrm{sub}(v)\backslash \{(x,t):x\in\partial\Omega, t<\psi\}$. By Remark~\ref{Density:T}, for any $a\in \mathrm{gr}(\psi)$, $\Theta(\|T_k\|, a)=\F{1}{2}$. Since $T$ is also a $\Lambda$-minimizing current, and since $\Sigma \backslash \mathrm{gr}(\psi)\cap B_r(a)$ has only one component for sufficiently small $r$, it follows from~\cite[Theorem~4.5(2)]{DS93} that
        \be \label{key:fact}
        \Theta(\|T\|, a)=\F{1}{2}.
        \ene

        We now proceed similarly to the proof of~\cite[Lemma~3.6]{zhou22a}. Fix a sufficiently small $r_1>0$. Note that all $\spt(T_k)\backslash \mathrm{gr}(\psi)$ are $C^2$ hypersurfaces with bounded mean curvature. By~\cite[Section~3.4]{All75}, there exists a positive constant $\kappa_1$ such that the function
        \be\label{eq:monotone}
        e^{\kappa_1 r} \F{\mathrm{H}^n(B_r(a)\cap \spt(T_k))}{\omega_nr^n}
        \ene
        is decreasing for any $r\in (0,r_1)$, any $k\geq 1$. By the codimension-one Allard regularity theorem with $C^{1,\alpha}$ boundary~\cite[p.\,418]{All75}, for any fixed $\gamma\in (0,\F{\alpha}{2})$ there exists $\delta>0$ such that if
        \be \label{det:A:B}
        \F{\mathbf{M}_{B_r(a)}(T_k)}{\omega_nr^n}=\F{\mathrm{H}^n(B_r(a)\cap \spt(T_k))}{\omega_nr^n}\leq \F{1}{2}+\delta,
        \ene
        then $\spt(T_k)\cap B_r(a)$ are uniformly $C^{1,\gamma}$ graphical submanifolds in $\R^{n_0}$ with boundary $\mathrm{gr}(\psi)$ which depend only on $\Lambda$, $\delta$ and $n$.

        Suppose, for contradiction, that no such positive constant $\delta$ exists for all $T_k, k\geq 1$. Then there is a $\delta_0>0$, a subsequence $\{k'\}_{k'=1}^\infty$ and radii $r_{k'}\rightarrow 0$ such that
        $$
        \F{\mathrm{H}^n(B_{r_{k'}}(a)\cap \spt(T_{k'}))}{\omega_nr_{k'}^n}\geq \F{1}{2}+\delta_0.
        $$
        By the monotonicity in~\eqref{eq:monotone}, there exists an $r_2<r_1$ independent of $k'$ such that for all $r\in (r_{k'}, r_2)$ and $k'\geq 1$ it holds that
        $$
        \F{\mathrm{H}^n(B_r(a)\cap \spt(T_{k'}))}{\omega_nr^n}\geq \F{1}{2}+\F{\delta_0}{2}.
        $$
        Arguing as in the proof of~\cite[Lemma~2.11(2)]{zhou22a}, there exists a dense set $I$ in $(0,r_2)$ such that for any $r\in I$,
        $$\lim_{k'\rightarrow \infty} \mathrm{H}^n(\partial B_r(a)\cap \spt(T_{k'})) =\mathrm{H}^n(\partial B_r(a)\cap\Sigma ) =0.$$
        By the convergence of BV functions (see, e.g.,~\cite{Giu84}), we have
        $$
        \F{\mathrm{H}^n(B_r(a)\cap \Sigma)}{\omega_nr^n}\geq \F{1}{2}+\F{\delta_0}{2}
        $$
        for all $r\in (0, r_2)$, which contradicts~\eqref{key:fact}. Therefore, there exist an $r>0$ and $\delta$ such that~\eqref{det:A:B} holds for all $k\geq 1$. For any $\gamma\in (0,\F{\alpha}{2})$, conclusion (1) follows from the Allard regularity theorem with $C^{1,\alpha}$ boundary~\cite[p.\,419]{All75}.

        To prove conclusion (2), by conclusion (1) and Lemma~\ref{lm:slice:convergence}, it suffices to show that $\Sigma\backslash \mathrm{gr}(\psi)$ is $C^{1,\gamma}$ and embedded for $2\leq n\leq 7$. By Theorem~\ref{thm:regularity}, if $2\leq n\leq 6$, this claim is true because $\mathrm{sub}(v)$ is a $\Lambda$-perimeter minimizer in $\Omega_{\delta}\PLH \R\backslash \mathrm{gr}(\psi)$ by item (1) of Theorem~\ref{thm:lambda:perimeter}.

        The remaining case is that $n=7$. The dimension of $N\PLH(a,b)$ is 8. By Theorem~\ref{thm:regularity}, the singular set of $\Sigma \backslash \mathrm{gr}(\psi)$ is discrete. Suppose it is nonempty, and let $p_0=(x_0, v_0)$ be a singular point. Then for any $\lambda>0$, define the map $T_\lambda: \R^{n_0+1}\rightarrow \R^{n_0+1}$ by $T_\lambda (p)=\F{p-p_0}{\lambda}$. Since $\mathrm{sub}(v)$ is a $\Lambda$-perimeter minimizer in a neighborhood of $p_0$, as $\lambda\rightarrow 0$, $T_\lambda (\mathrm{sub}(v))$ will converge to an $n$-dimensional minimizing cone $C$ in $T_{p_0}(N\PLH(a,b))\cong \R^{n+1}$. Note that
        $$
        \begin{aligned}
            T_\lambda( \mathrm{sub}(v)) &=\{(y,r):y \in \F{N-x_0}{\lambda},  r<v_\lambda(y) \},\\ v_\lambda(y) &= \F{v(y\lambda+x_0)-v(x_0)}{\lambda}.
        \end{aligned}
        $$

        As $\lambda\rightarrow 0$, $\F{N-x_0}{\lambda}$ converges locally smoothly to $T_{x_0}N$, a linear $n$-dimensional space of $T_{p_0}(N\PLH(a,b))$. Since $0\in \partial C$ and $n=7$, the point $0$ is the only possible singular point of $\partial C$. The normal vector of $T_\lambda(\partial \mathrm{sub}(v))$ converges pointwise to $\vec{\nu}_c$, the normal of $\partial C\setminus \{0\}$. Let $\partial_r$ be the coordinate vector field of $\R$. On each $T_{\lambda}(\partial \mathrm{sub}(v)\setminus \{p_0\})$, $\langle \vec{\nu},\partial_r\rangle \geq 0$, and so $\langle \vec{\nu}_c, \partial_r \rangle \geq 0$ on $\partial C\setminus \{0\}$. Since $\partial C\setminus \{0\}$ is smooth and minimal, it follows that $\langle \vec{\nu}_c,\partial_r\rangle$ satisfies
        $$
        \Delta \langle \vec{\nu}_c, \partial_r \rangle +|A|^2\langle \vec{\nu}_c, \partial_r\rangle \equiv 0
        $$
        on $\partial C\setminus \{0\}$, where $|A|^2$ is the norm of the second fundamental form and $\Delta$ is the Laplacian operator of $\partial C\setminus \{0\}$. For a derivation, see~\cite[Lemma~2.2]{zhou19a}. By the Harnack inequality, only two cases happen: $\langle \vec{\nu}_c, \partial_r\rangle\equiv 0$ or $\langle \vec{\nu}_c, \partial_r\rangle>0$ in $\partial C\setminus \{0 \}$ in $\R^{n+1}$. The first case is impossible because $0$ is an isolated singular point of $\partial C$. In the second case, $\partial C$ is a minimal graph with a possible singular point. By~\cite[Theorem~1]{Simon77}, such an isolated point is removable and $\partial C$ is smooth. Because $\partial C$ is a smooth cone passing through $0$, $C$ is contained in a half-space. By~\cite[Theorem~15.5]{Giu84}, $C$ is a half-space in $T_{p_0}(N\PLH (a,b))$. By the Allard regularity theorem, $\partial \mathrm{sub}(v)$ is $C^{1,\gamma}$ near $p_0$ for any $\gamma\in (0,\F{\alpha}{2})$. As a result, $\Sigma\setminus \mathrm{gr}(\psi)$ is $C^{1,\gamma}$ and embedded for $n=7$.

        This completes the proof of conclusion (2).
    \end{proof}

    \section{Existence of PMC hypersurfaces}\label{sec:A}
    In this section, we generalize Gerhardt's result to the case of conformal product manifolds with fixed (possibly empty) boundaries. As discussed in Theorem~\ref{at:G:key}, such existence is closely related to solving the following PMC equation
    \be \label{def:L:operator}
    \mathcal{L}_\mathcal{H}(u): =-\div\Bigl(\F{Du}{\omega}\Bigr)-\mathcal{H}\Bigl(x,u,-\F{Du}{\omega},\F{1}{\omega}\Bigr)
    \ene
    where $\omega=\sqrt{1+|Du|^2}$, $\mathcal{H}(p ,X)=\mathcal{H}(x,r, Y,t)$ is a $C^{1,\alpha}$ function on the tangent bundle $T (N\PLH(a,b))$ for $p=(x,r)\in N\PLH\R$ and $X=Y+t\partial_r\in T(N\PLH\R)$ with $x\in N,r,t\in \R$ and $Y\in T_x N$.

    The main result of this section is formulated as follows.

    \begin{theorem}\label{thm:main}
        Let $N$ be an $n$-dimensional closed or $C^{2,\alpha}$ compact Riemannian manifold with interior $N^0$ and metric $\sigma$, where $2\leq n\leq 7$.

        Let $M_f$ be the conformal product manifold defined as in~\eqref{def:A:B}. Let $\mathcal{H}$ be a $C^{1,\alpha}$ function on $TM_f$, and assume there exist two functions $u_1, u_0$ in $C^{2,\alpha}(N^0)\cap C(N)$ such that $u_1<u_0$ on $N^0$, and the pair $(\mathrm{gr}(u_1),\mathrm{gr}(u_0))$ forms a barrier for $(M_f, \mathcal{H})$ in the sense of Definition~\ref{def:A}, i.e.,
        $$
        \rH_{\mathrm{gr}(u_1)}-\mathcal{H}\leq 0,\quad  \rH_{\mathrm{gr}(u_0)}-\mathcal{H}\geq 0,
        $$
        where $\rH_{\mathrm{gr}(u_i)}$ denotes the mean curvature of $\mathrm{gr}(u_i)$ with respect to the upward unit normal vector for $i=0,1$. If $\partial N\neq \emptyset$, we require that $u_1=u_0=\psi$ on $\partial N$ for some $\psi \in C^{1,\alpha}(\partial N)$.

        Then there exists a $C^{3,\alpha}$ embedded orientable hypersurface $\Sigma\subset M_f$, homeomorphic to $N$, with mean curvature $\mathcal{H}$. If $\partial N\neq \emptyset$, then $\partial \Sigma=\mathrm{gr}(\psi)$ is of class $C^{1,\beta}$. If $N$ is closed, then $\Sigma$ is also closed. Here $\beta$ is any number in $(0,\F{\alpha}{2})$.
    \end{theorem}

    \begin{Rem}\label{rk:description}
        In the proof, we shall show that there exists a $v\in W^{1,1}(N^0)$ with $u_1\leq v\leq u_0$, $v\in \mathcal{D}_{\kappa}(N^0)$ for some positive constant $\kappa$ from Theorem~\ref{key:thm}, and that $\Sigma$ is given by $\partial \mathrm{sub}(v)\cap (N^0\PLH \R)$. However, $\Sigma$ may not be a $C^1$ graph over $N^0$.
    \end{Rem}

    \begin{proof}
        By Theorem~\ref{at:G:key}, finding a $C^{3,\alpha}$ embedded hypersurface in $M_f$ with mean curvature $\mathcal{H}$ is equivalent to finding a $C^{3,\alpha}$ embedded hypersurface in the product manifold $N\PLH (a,b)$ with mean curvature $\mathcal{H}'(p,\vec{\nu})$, where $\mathcal{H}':T(N\PLH (a,b))\rightarrow \R$ is defined by
        $$
        \mathcal{H}'(x,r,Y,t)=e^{f(x,r)}\{\mathcal{H}(x,r,Y,t)-n \langle \nabla f, Y\rangle-n\F{\partial f}{\partial r}t\}.
        $$
        Here, $p=(x,r)$, $X=(Y+t\partial_r )$, $Y\in T_x N$, and $\partial_r$ is the coordinate vector field on $(a,b)$. The condition that $(\mathrm{gr}(u_1),\mathrm{gr}(u_0))$ forms a barrier for $(M_f, \mathcal{H})$ is equivalent to $(\mathrm{gr}(u_1), \mathrm{gr}(u_0))$ forming a barrier for $(N\PLH (a,b), \mathcal{H}')$, i.e.,
        \be \label{var:barrier:condition}
        \mathcal{L}_{\mathcal{H}'}(u_1)\leq  0, \quad  \mathcal{L}_{\mathcal{H}'}(u_0)\geq 0, \quad  u_1<u_0 \text{ on } N^0,
        \ene
        where $\mathcal{L}_{\mathcal{H'}}(u)$ is given by~\eqref{def:L:operator}. Therefore, it suffices to show Theorem~\ref{thm:main} with $\mathcal{H} \equiv \mathcal{H}'$ in the product manifold ($f\equiv 0$) under the barrier condition in~\eqref{var:barrier:condition}.

        From now on assume $f\equiv 0$, $\mathcal{H}\equiv \mathcal{H}'$, and the barrier condition~\eqref{var:barrier:condition} holds. Note that there exists an interval $(c_1, c_2)\subset \subset (a,b)$ satisfying $c_1\leq u_1\leq u_0\leq c_2$. Let $h(r):(a,b)\rightarrow [0,1]$ be a smooth function satisfying $h(r)\equiv 1$ for any $r\in [c_1, c_2]$, with $\spt(h)\subset \subset (a,b)$ and $0\leq h\leq 1$.

        We define a new operator $\mathcal{L}_{2}(u)$ by
        \be \label{eq:one:p}
        \mathcal{L}_{2}(u):=-\div\Bigl(\F{Du}{\omega}\Bigr)-h(u)\mathcal{H}\Bigl(x,u,-\F{Du}{\omega},\F{1}{\omega}\Bigr)+\gamma(u-u_1),
        \ene
        where $\omega=\sqrt{1+|Du|^2}$ and $\gamma$ is a constant satisfying
        \be \label{maximum:condition}
        -\F{\partial (h \mathcal{H})}{\partial r}(x,r, X, t)+\gamma\geq 1.
        \ene
        By the definition of $h(r)$, it follows that on $N^0$,
        $$
        \mathcal{L}_{2}(u_1)=\mathcal{L}_{\mathcal{H}'}(u_1)\leq 0,\quad  \mathcal{L}_{2}(u_0)\geq \mathcal{L}_{\mathcal{H}'}(u_0)\geq 0.
        $$
        Let $\mathcal{F}_2(x,r,X,t)=h(r)\mathcal{H}(x,r,X,t)-\gamma(r-u_1)$. Then $(\mathrm{gr}(u_1), \mathrm{gr}(u_0))$ forms a barrier for $(N\PLH(a,b),\mathcal{F}_2)$ with $\F{\partial \mathcal{F}_2}{\partial r}\leq 0$. By Theorem~\ref{thm:PDE:existence:}, there exists a function $u_2\in C^{2,\alpha}(N^0)\cap C(N)$ such that $\mathcal{L}_2(u)\equiv 0$ with $u_1\leq u_2\leq u_0$.

        From now on, fix $\gamma$ in~\eqref{maximum:condition}. By induction on $m$, we will define a $C^{2,\alpha}$ sequence $\{u_m\}_{m=2}^\infty$ satisfying
        \be \label{condition:DS}
        u_{m-1}\leq u_m\leq u_0,\quad \mathcal{L}_m(u_m)=0 \quad \text{on } N^0.
        \ene
        Suppose that the operator $\mathcal{L}_{m-1}$, function $\mathcal{F}_{m-1}$ and $u_{m-1}$ have been defined for $m \ge 3$ to satisfy
        $$
        \frac{\partial \mathcal{F}_{m-1}}{\partial r} \le 0, \quad \mathcal{L}_{m-1}(u_{m-1})=0, \quad u_{m-2}\le u_{m-1} \le u_0.
        $$
        Let $\mathcal{F}_m:=h \mathcal{H}-\gamma(r-u_{m-1})$ and define
        \begin{equation} \label{def:L_k}
            \mathcal{L}_m(u):=-\div\Bigl(\F{Du}{\omega}\Bigr)-\mathcal{F}_m\Bigl(x,u,-\F{Du}{\omega},\F{1}{\omega}\Bigr).
        \end{equation}
        Since $\mathcal{F}_{m}=\mathcal{F}_{m-1}+\gamma(u_{m-2}-u_{m-1})$, it holds that $\frac{\partial \mathcal{F}_m}{\partial r} \le 0$, and
        \begin{equation}
            \mathcal{L}_m(u)=\mathcal{L}_{m-1}(u)-\gamma(u_{m-1}-u_{m-2}).
        \end{equation}
        Then $\mathcal{L}_m(u_{m-1})\le 0$, $\mathcal{L}_{m}(u_0) \ge 0$. Therefore, $(\mathrm{gr}(u_{m-1}), \mathrm{gr}(u_0))$ forms a barrier for $(N \PLH \R, \mathcal{F}_m)$. By Theorem~\ref{thm:PDE:existence:}, there exists a function $u_m \in C^{2,\alpha}(N^0)\cap C(N)$ to satisfy $\mathcal{L}_m(u_m)=0$ and $u_{m-1} \le u_m \le u_0$ on $N^0$ for any $m \ge 1$. The sequence $\{u_m\}_{m=0}^\infty$ is thus well-defined.

        By~\eqref{condition:DS}, there exists a positive constant $\kappa$, depending on $u_1$, $u_0$ and $\gamma$, such that
        \be \label{fact:uk}
        \max\{|u|, |\div(\F{Du}{\sqrt{1+|\nabla u|^2}})|\}\leq \kappa \text{ on } N^0
        \ene
        for all $u=u_m$, $m=2,3,\cdots$.

        By~\eqref{fact:uk} and Theorem~\ref{key:thm}, the sequence $\{u_m\}^\infty_{m=2}$ converges to a function $v \in W^{1,1}(N^0)$ with $u_1\leq v \leq u_0$, and $v$ is $C^{1,\gamma}$ in an open dense set $N_v \subset N^0$ for any $\gamma\in (0,\F{1}{2})$. Now define the hypersurface $\Sigma$ by
        \be \label{def:sigma}
        \Sigma :=\partial \mathrm{sub}(v)\cap(N^0\PLH \R)\subset \{(x,t):x\in N^0, u_1(x)<t<u_0(x)\}.
        \ene

        Combining the fact that $2\leq n \leq 7$ with Theorem~\ref{thm:boundary:convergence} and~\eqref{fact:uk}, $\Sigma$ is a $C^{1,\alpha}$ embedded hypersurface with (possibly empty) $C^{1,\beta}$ boundary $\mathrm{gr}(\psi)$ for any $\beta \in (0,\F{\alpha}{2})$. Note that item (2) of Theorem~\ref{thm:boundary:convergence} still holds even when $N$ is closed. Therefore, $\mathrm{gr}(u_m)$ converges uniformly to $\Sigma$ in the sense of the $C^{1,\beta}$ norm. As a result, $\Sigma$ is homeomorphic to $N$, whether $N$ is closed or a $C^{2,\alpha}$ compact manifold.

        It remains to compute the mean curvature of $\Sigma$. We first consider the portion $(N_v\PLH \R)\cap \Sigma$, which is a $C^{1,\gamma}$ graph over $N_v$. For any $x\in N_v$, $\nabla v$ is finite and continuous in a neighborhood of $x$. By item (1) of Lemma~\ref{lm:slice:convergence}, the graph of $u_m$ converges uniformly to $\Sigma$ in the $C^{1,\gamma}$ norm on a neighborhood of $(x,v(x))$ for some $\gamma \in (0,\F{1}{2})$. Choose an embedded ball $B_r(x)\subset N$ centered at $x$ with small radius $r$. Due to the $C^{1,\gamma}$ convergence of $\mathrm{gr}(u_k)$, $\{\nabla u_k\}_{k=1}^\infty$ converge uniformly to $\nabla v$ in $B_r(x)$. As a result, $\max_{B_{\F{r}{2}}(x) }|\nabla u_k|$ is uniformly bounded. From~\eqref{condition:DS} and~\eqref{def:L_k}, the classical Schauder estimates imply that their $C^2$ norms are uniformly bounded on $B_{\F{r}{2}}(x)$, i.e.,
        \be
        \|u_m\|_{C^2(B_{\F{r}{2}}(x))}\leq C, \quad m\geq 2,
        \ene
        where $C$ is a constant independent of $m$. As a result, $\Sigma\cap (B_{\F{r}{2}}(x)\PLH \R)$ is a $C^2$ graph over $B_{\F{r}{2}}(x)$. Moreover, the mean curvature of $\Sigma \cap B_{\F{r}{2}}(x)$ with respect to the upward normal vector is given by:
        \be
        \begin{split}
            \rH_{\Sigma}((y,v(y)))&=-\div\Bigl(\F{\nabla v}{\omega}\Bigr)=\lim_{m\rightarrow \infty }\mathcal{H}\Bigl(y,u_m, \F{\nabla u_m}{\omega},\F{1}{\omega}\Bigr)+\gamma\F{(u_m-u_{m-1})}{\omega}\\
            &=\mathcal{H}\Bigl(y,v(y),-\F{Dv}{\omega},\F{1}{\omega}\Bigr),
        \end{split}
        \ene
        for any $y\in B_{\F{r}{2}}(x)$, where $\omega=\sqrt{1+|Dv|^2}$ (and analogously for $u_m$). Here we use the fact $\lim\limits_{m\to \infty}(u_m-u_{m-1})=0$ on $B_{\F{r}{2}}(x)$, which follows from the $C^{1,\gamma}$ convergence of the graph $\mathrm{gr}(u_m)$. Since $x$ is arbitrary, we conclude that
        \be \label{det:dense:mc}
        \rH_{\Sigma}(p)=\mathcal{H} (p, \vec{\nu})
        \ene
        for any $p\in \Sigma\cap (N_v\PLH \R)$, where $\vec{\nu}$ is the normal vector of $\Sigma$.
        Since $\mathcal{H}$ is $C^{1,\alpha}$, it follows that $v$ is $C^{3,\alpha}$ over $N_v$, and hence $\Sigma\cap ( N_v\PLH \R)$ is also $C^{3,\alpha}$.

        Next we compute the mean curvature of $\Sigma\setminus (N_v\PLH \R)$, where its normal vector $\vec{\nu}$ satisfies $\langle \vec{\nu},\partial_r\rangle \equiv 0$.
        Fix any point $p \in \Sigma\setminus (N_v \PLH \R)$. We choose a new coordinate chart
        \be
        \{W,x=(x_1,x_2,\cdots, x_n,x_{n+1})\},
        \ene
        in a neighborhood of $p$ such that $\partial_{n+1}(p)=\vec{\nu}(p)$. Then $\Sigma\cap W$ can be expressed as the graph of a $C^{1,\gamma}$ function $f(x_1, \cdots, x_n)$,
        \begin{equation}
            \Sigma\cap W=\{(p',x_{n+1}): x_{n+1}=f(p')\text{ with } p'=(x_1,\cdots, x_n) \in V\}
        \end{equation}
        where $V\subset \R^n$ is a sufficiently small neighborhood of the origin $0$ in $\mathbb{R}^n$ and $p=(0, f(0))$. With this new coordinate chart, the graph of $f$ is the limit of a sequence of graphs of ${f_k}$, which are also the graphs of the functions $u_m$ in the original product metric coordinate $\{(x,t):x\in N,t\in \R\}$. Moreover, the mean curvature of these graphs is uniformly bounded.

        By Theorem~\ref{key:thm} and Lemma~\ref{lm:fullmeasure}, $f$ belongs to $W^{1,1}(V)$, and $\Sigma\cap (N_v\PLH \R)$ is open dense in $\Sigma$ with respect to the $n$-dimensional Hausdorff measure. As a result, $f$ is $C^{3,\alpha}$ except for a set of $n$-dimensional Hausdorff measure zero in $V$, while $f\in C^{1,\gamma}(V)$. Equation~\eqref{det:dense:mc} is equivalent to
        \be \label{det:mc:sense:A}
        \partial_i\bigl(a^{ij}(x,f, \nabla f)f_{i}\bigr)+b^i(x, f, \nabla f)f_j+\mathcal{H}''(x,f, \nabla f)=0,
        \ene
        which is a second order elliptic equation of divergence type with smooth coefficients, and $\mathcal{H}''$ is a $C^{1,\alpha}$ function derived from $\mathcal{H}$. Moreover, $f\in W^{1,1}(V)\cap C^{1,\gamma}(V)$ is a weak solution to equation~\eqref{det:mc:sense:A}. By~\cite[Theorem~10.1]{LU68}, $f\in W^{2,2}(V)$, and then~\cite[Theorem~12.1]{LU68} implies that $f\in C^{2,\alpha}(V)$. Hence, $\Sigma$ is a $C^{2,\alpha}$ embedded hypersurface in $W$, and the mean curvature of $\Sigma$ is continuous with respect to $p$ and $\vec{\nu}$ in $W$. Since $\Sigma\cap (N_v\PLH \R)$ is dense in $\Sigma$,~\eqref{det:dense:mc} holds throughout $\Sigma\cap W$.

        Since $p$ is arbitrarily chosen in $\Sigma\setminus (N_v\PLH\R)$, the formula
        \be \label{mc:equation}
        \rH_{\Sigma}(p)=\mathcal{H}(p, \vec{\nu})
        \ene
        holds for any point $p$ in $\Sigma$. This completes the proof.
    \end{proof}

    \section{Quasi-decreasing case}\label{sec:B}
    The PMC hypersurface $\Sigma$ in Theorem~\ref{thm:main} may not be a $C^1$ graph over $N$. In this section, we introduce a condition on the PMC function $\mathcal{H}$---distinct from the decreasing condition in~\eqref{condition:monotone}---to ensure that $\Sigma$ is a $C^1$ graph.

    \begin{Def}
        A function $\mathcal{H}(x,r,X,t):T(N\PLH \R)\rightarrow \R$ is called quasi-decreasing if it can be expressed as
        \be \label{quasi:monotone}
        \mathcal{H}(x,r,X,t)= \mathcal{H}_1(x,r, X,t)+t\mathcal{H}_2(x,r,X,t), \quad \F{\partial \mathcal{H}_1}{\partial r}\leq 0,
        \ene
        where $\mathcal{H}_1, \mathcal{H}_2$ are $C^{1,\alpha}$ functions on the tangent bundle $T (N\PLH(a,b))$.
    \end{Def}

    The operator~\eqref{def:L:operator} then takes the form
    \be \label{def:L:operator:A}
    \mathcal{L}_{\mathcal{H}}(u): =-\div\Bigl(\F{Du}{\omega}\Bigr)-\mathcal{H}_1\Bigl(x,u,-\F{Du}{\omega},\F{1}{\omega}\Bigr)-\F{1}{\omega}\mathcal{H}_2\Bigl(x,u,-\F{Du}{\omega},\F{1}{\omega}\Bigr),
    \ene
    where $\omega=\sqrt{1+|Du|^2}$.

    \begin{theorem}\label{quasi:monotone:case}
        Let $\mathcal{L}_{\mathcal{H}}$ be the operator defined in~\eqref{def:L:operator:A}, where $\mathcal{H}$ is quasi-decreasing as in Definition~\ref{quasi:monotone}. Assume the following:
        \begin{enumerate}
            \item $N$ is an $n$-dimensional ($7\geq n\geq 2$) closed or $C^{2,\alpha}$ compact Riemannian manifold with interior $N^0$;
            \item There exist two functions in $C^2(N^0)\cap C(N)$, $u_1<u_0$, such that $(\mathrm{gr}(u_1), \mathrm{gr}(u_0))$ forms a barrier for $(N\PLH\R, \mathcal{H})$ as in Definition~\ref{def:A}. Here $N\PLH \R$ is the product manifold.
        \end{enumerate}
        Then there is a function $v\in C^{3,\alpha}(N^0)\cap C(N)$ satisfying $u_1 \leq v\leq u_0$ and
        $$
        \mathcal{L}_{\mathcal{H}}(v)\equiv 0 \text{ on } N^0.
        $$
    \end{theorem}

    \begin{proof}
        By Theorem~\ref{thm:main}, there exists a $C^{3,\alpha}$ embedded hypersurface $\Sigma$ with mean curvature $\mathcal{H}$. By Remark~\ref{rk:description}, $\Sigma=\partial \mathrm{sub}(v)\cap (N^0\PLH \R)$ for $v\in W^{1,1}(N)$ with $u_1\leq v\leq u_0$.

        To show the conclusion, it is sufficient to show that $\langle \vec{\nu},\partial_r\rangle >0$ on $\Sigma$, where $\vec{\nu}$ is the upward normal vector of $\Sigma$. This implies that $v \in C^1(N^0)$. Since $\mathcal{H}$ is $C^{1,\alpha}$, it then follows that $v\in C^{3,\alpha}$ and satisfies $\mathcal{L}_{\mathcal{H}}(v)=0$.

        Let $\langle\cdot,\cdot\rangle$ be the inner product of the product manifold $N\PLH\R$ and define $\Theta =\langle \vec{\nu}, \partial_r \rangle$, which is non-negative on $\Sigma$. As in~\cite[Lemma~2.2]{zhou19a}, the following identity holds
        \be \label{de:A:hold}
        \Delta \Theta+(|A|^2+\bar{R}ic(\vec{\nu},\vec{\nu}))\Theta-\langle \nabla H_{\Sigma}, \partial_r\rangle=0
        \ene
        where $\Delta, \nabla$ are the Laplacian and the covariant derivative of $\Sigma$, respectively, $\vec{\nu}$ is the upward unit normal vector of $\Sigma$, $\rH_{\Sigma}$ is the mean curvature of $\Sigma$ with respect to $\vec{\nu}$, and $\bar{R}ic$ is the Ricci curvature of the product manifold $N\PLH \R$.

        Let $\partial_r^T$ be the tangential component of $\partial_r$ in $T\Sigma$. By definition
        \be
        \partial_r^T:=\partial_r-\Theta\vec{\nu}.
        \ene
        By Theorem~\ref{key:thm}, $v$ is $C^{1,\gamma}$ on an open dense set $N_v\subset N^0$. The identity $\rH_{\Sigma}=\mathcal{H}(p,\vec{\nu})$ implies that $v$ is at least a $C^3$ function on $N_v$. Thus, on $\Sigma\cap (N_v\PLH\R)$, $\rH_{\Sigma}=\mathcal{H}(x,v,-\F{\nabla v}{\omega},\F{1}{\omega})$ is a $C^2$ function on $N_v$. As a result, on $\Sigma\cap(N_v\PLH\R)$, we obtain:
        \be\label{formula:theta}
        \Delta \Theta+(|A|^2+\bar{R}ic(\vec{\nu},\vec{\nu})+\langle \nabla v, \nabla\mathcal{H}\rangle \Theta)\Theta=0,
        \ene
        where $\vec{\nu}=\F{-\nabla v+\partial_r}{\omega}$ and $\omega=\sqrt{1+|\nabla v|^2}$. By definition, $\Theta =\F{1}{\omega}$.

        We now follow the proof of~\cite[Theorem~A.3]{zhou22a}. Let $\{\partial_1, \cdots, \partial_n \}$ be a local frame on $N_v$, $\sigma_{ij}=\langle \partial_i, \partial_j\rangle$, $(\sigma^{ij})=(\sigma_{ij})^{-1}$. Let $v_{i}$ and $v_{ij}$ denote the first and second covariant derivatives of $v$, respectively. Define $v^k=\sigma^{ki}v_i$. The gradient of $v$ is given by $\nabla v=v^k\partial_k$. Then $\tilde{v}^k=\Theta v^k$, and $\vec{\nu}=\Theta\partial_r-\tilde{v}^k\partial_k$, where $\Theta=\F{1}{\sqrt{1+|\nabla v|^2}}$. The set $\{X_i =\partial_i +v_i \partial_r,i=1,\cdots, n\}$ forms a frame for $T\Sigma$. The metric of $\Sigma$ is $(g_{ij})=(\langle X_i, X_j\rangle) =(\sigma_{ij}+v_iv_j)$, with its inverse matrix $(g^{ij})=(\sigma^{ij}-\F{v^iv^j}{\omega^2})$.

        The following two formulas are useful (see~\cite[section~2]{zhou18}).
        \begin{align*}
            \tilde{v}^k_i=\F{1}{\omega}g^{kl}v_{li}, \quad |A|^2=\tilde{v}^k_i\tilde{v}_k^i.
        \end{align*}
        Let $Y=p^1\partial_1+\cdots+p^n\partial_n$. We express $\mathcal{H}(x,r,Y,t)$ as $\mathcal{H}(x,r, p^1,\cdots, p^n, t)$, $x=(x^1,\cdots, x^n)$ and $\partial_k$ denotes $\F{\partial}{\partial x_k}$. Note that on $\Sigma\cap (N_v\PLH\R)$, it holds that
        \be
        \mathcal{H}|_{\Sigma}=\mathcal{H}\Bigl(x,v,-\frac{\nabla v}{\omega},\frac{1}{\omega}\Bigr)=\mathcal{H}_1+\frac{1}{\omega}\mathcal{H}_2
        \ene
        with
        \be
        \frac{\partial \mathcal{H}}{\partial t}=\frac{\partial \mathcal{H}_1}{\partial t}+t\frac{\partial \mathcal{H}_2}{\partial t}+\mathcal{H}_2,\quad \frac{\partial \mathcal{H}}{\partial r}=\frac{\partial \mathcal{H}_1}{\partial r}+t\frac{\partial \mathcal{H}_2}{\partial r}.
        \ene
        Then:
        \be
        \begin{split}
            -\Theta^2 \langle \nabla \mathcal{H}, \nabla v\rangle &=-\Theta \tilde{ v}^k \partial_k (\mathcal{H}) \\
            &=-\Theta \tilde{v}^k\Bigl(\F{\partial \mathcal{H}}{\partial x_k}+\F{\partial \mathcal{H}}{\partial r} v_k+\F{\partial \mathcal{H}}{\partial p_{j}}\tilde{v}^j_k+\F{\partial \mathcal{H}}{\partial t}\Theta_k\Bigr)\\
            &=-\Theta \Bigl( \tilde{v}^k\F{\partial \mathcal{H}}{\partial x_k}+\F{\partial \mathcal{H}}{\partial r} \tilde{v}^kv_k+\F{\partial \mathcal{H}}{\partial p_{j}}\tilde{v}^j_k\tilde{v}^k\Bigr)-\F{\partial \mathcal{H}}{\partial t}\langle \nabla\Theta, \partial_r \rangle \notag\\
            &=-\Theta \Bigl(\tilde{v}^k\F{\partial \mathcal{H}}{\partial x_k}+\F{\partial \mathcal{H}}{\partial p_j}\tilde{v}^j_k\tilde{v}^k\Bigr) -\Theta^2|Dv|^2\F{\partial \mathcal{H}_1}{\partial r}-\Theta\F{|Dv|^2}{1+|Dv|^2}\F{\partial \mathcal{H}_2}{\partial r}\\
            &\quad-\F{\partial \mathcal{H}}{\partial t}\langle \nabla\Theta, \partial_r \rangle.
        \end{split}
        \ene
        on $N_v$.

        Because $v$ satisfies $a\leq u_1\leq v\leq u_0\leq b$ for some interval $[a,b]$, there is a positive number $\kappa$ only depending on $\mathcal{H}_1, \mathcal{H}_2$ and $N\PLH [a,b]$ such that on $\Sigma$
        \be
        \max\{|\tilde{v}^k\F{\partial \mathcal{H}}{\partial x_k}|, |\F{\partial \mathcal{H}_1}{\partial r}|,|\F{\partial \mathcal{H}_2}{\partial r}|, |\F{\partial \mathcal{H}}{\partial p_i}|, |\F{\partial \mathcal{H}}{\partial t}|, |\bar{R}ic_N(\vec{\nu},\vec{\nu})|,\Theta \}\leq \kappa.
        \ene
        Note that $\kappa$ is independent of the set $N_v$.

        As a result, for any fixed $\varepsilon>0$ we have
        $$
        \F{\partial \mathcal{H}}{\partial p_i}\tilde{v}^i_k\tilde{v}^k\geq -\varepsilon |A|^2 -\F{\kappa}{4\varepsilon}.
        $$
        By~\eqref{quasi:monotone}, $\F{\partial \mathcal{H}_1}{\partial r}\leq 0$. Therefore,
        \be
        -\Theta^2 \langle D \mathcal{H}, Dv\rangle\geq -\Theta \Bigl(\kappa\bigl(2+\F{1}{4\varepsilon}\bigr) +\varepsilon|A|^2\Bigr)-\F{\partial \mathcal{H}}{\partial t}\langle \nabla\Theta,\partial_r\rangle.
        \ene
        This is the only place that the quasi-decreasing condition in~\eqref{quasi:monotone} was used. Substituting this into~\eqref{formula:theta} gives that on $\Sigma\cap (N_v\PLH\R)$

        \be \label{det:atD}
        \Delta \Theta+\bigl(|A|^2(1-\varepsilon)-\kappa\bigl(3+\F{1}{4\varepsilon}\bigr)\bigr)\Theta-\bigl(\mathcal{H}_2+\F{\partial \mathcal{H}_1}{\partial t}\bigr)\langle \nabla \Theta, \partial_r\rangle \leq 0.
        \ene

        Since $\Sigma$ is $C^{3,\alpha}$, $\Theta$ is a $C^{2,\alpha}$ function on $\Sigma$. By Theorem~\ref{key:thm}, $\Sigma\cap (N_v\PLH\R)$ is dense in $\Sigma$, and $\kappa$ and $\varepsilon$ are constants independent of $N_v$. Then the inequality~\eqref{det:atD} holds on all of $\Sigma$. Because $\Theta\geq 0$ on $\Sigma$, by the Harnack inequality, either $\Theta >0$ or $\Theta\equiv 0$ on $\Sigma$. Because $N_v$ is open dense, only the former case can happen. This means that $v\in C^{3,\alpha}(N^0)$. And the graph of $v$ coincides with $\Sigma$, satisfying $\mathcal{L}_{\mathcal{H}}(v)\equiv 0$ on $N^0$.

        This completes the proof.
    \end{proof}

\end{document}